\documentclass{amsart}
\usepackage[utf8]{inputenc}
\usepackage{amsmath}
\usepackage{hyperref}
\usepackage[all]{xy}
\usepackage{amsfonts}
\usepackage{amssymb}
\usepackage{amsthm}
\usepackage{graphicx}
\usepackage{color}
\newtheorem{theorem}{Theorem}[section]
\newtheorem{corolario}[theorem]{Corollary}
\newtheorem{lema}[theorem]{Lemma}
\newtheorem{prop}[theorem]{Proposition}
\newtheorem{proposition}[theorem]{Proposition}
\newtheorem{defn}[theorem]{Definition}

\theoremstyle{definition}
\newtheorem{remark}[theorem]{Remark}

\newtheorem{example}[theorem]{Example}

\title{Analytic Semiroots for Plane Branches and Singular Foliations}
\author{Felipe Cano}
\address{Felipe Cano. Departamento de Álgebra, Análisis Matemático,
Geometr\'\i a y Topolog\'\i a. Universidad de Valladolid.
Paseo de Bel\'en 7,
47011 -- Valladolid, SPAIN}
\email{fcano@uva.es}
\author{Nuria Corral}
\address{Nuria Corral. Departamento de Matemáticas, Estadística y Computación. Universidad de Cantabria. Avda. de los Castros s/n, 39005 -- Santander, SPAIN}
\email{nuria.corral@unican.es}
\author{David Senovilla-Sanz}
\address{David Senovilla-Sanz. Departamento de Matemáticas, Estadística y Computación. Universidad de Cantabria. Avda. de los Castros s/n, 39005 -- Santander, SPAIN}
\email{david.senovilla@unican.es}
\thanks{The authors are supported by the Spanish research project  PID2019-105621GB-I00/AEI/10.13039/501100011033 funded by the Agencia Estatal de Investigación - Ministerio de Ciencia e Innovación. The third author is also supported by a predoctoral contract “Concepción Arenal” of the  Universidad de Cantabria}
\begin{document}
\begin{abstract} The analytic moduli of equisingular plane  branches has the semimodule of differential values as the most relevant system of discrete invariants. Focusing in the case of cusps, the minimal system of generators of this semimodule is reached by the differential values attached to the differential $1$-forms of the so-called standard bases. We can complete a standard basis to an extended one by adding a last differential $1$-form that has the considered cusp as invariant branch and the ``correct'' divisorial order. The elements of such extended standard bases have the ``cuspidal'' divisor as a ``totally dicritical divisor'' and hence they define packages of plane branches that are equisingular to the initial one. These are the analytic semiroots. In this paper we prove that the extended standard bases are well structured from this geometrical and foliated viewpoint, in the sense that the semimodules of differential values of the branches in the dicritical packages are described just by a truncation of the list of generators of the initial semimodule at the corresponding differential value. In particular they have all the same semimodule of differential values.
\end{abstract}
\begin{center}
\noindent{\bf This paper will appear in \\ Bulletin of the Brazilian Mathematical  Society \\
\url{https://doi.org/10.1007/s00574-023-00344-w}} \\
\vspace*{\baselineskip}
\end{center}
\maketitle
\tableofcontents
\section{Introduction}
The analytic classification of plane branches starts with Zariski \cite{zariski}, who pointed the importance of the differential values in this problem. The semimodule of differential values was extensively described by Delorme \cite{delorme}, although the complete analytic classification is due to A. Hefez and  M. E. Hernandes \cite{hefez2}.

Geometrically, the ``most interesting'' differential values are viewed as the contact $\nu_{\mathcal C}(\omega)$ of a given branch $\mathcal C $ with the foliations  defined by diffe\-rential $1$-forms $\omega$ without common factors in the coefficients. From the moduli view-point, the semimodule of differential values $\Lambda$ is interpreted as the ``discrete structure'' supporting the continuous part of the moduli. More precisely, the semimodule $\Lambda$ has a well defined basis $\{\lambda_j\}_{j=-1}^s$; so, it is reasonable to fix our attention in the differential forms that produce precisely the elements of the basis as differential values: these are the elements of the standard bases (for more details, see \cite{hefez1, hefez3}).

In this paper we focus in the case of cusps, that is, branches with a single Puiseux pair $(n,m)$. Our objective is to describe the cusps close to a cusp $\mathcal C$, in terms of a given standard basis $\mathcal H$ and the dicritical foliated behaviour of the elements of $\mathcal H$ in the final divisor $E$ of the reduction of singularities of $\mathcal C$. Let us precise this.

We consider a cusp $\mathcal C$ with Puiseux pair $(n,m)$. In view of Zariski Equisingularity Theory, we know that the semigroup $\Gamma=n\mathbb Z_{\geq 0}+m\mathbb Z_{\geq 0}$ of $\mathcal C$ is an equivalent data of the equisingularity class of $\mathcal C$. The differential values define a semimodule $\Lambda^\mathcal C$ over $\Gamma$, that will have a strictly increasing basis
$$
\lambda_{-1}=n,\lambda_0=m,\lambda_1,\ldots,\lambda_s,
$$
to be the minimal one such that $\Lambda^\mathcal C=\cup_{j=-1}^s(\lambda_j+\Gamma)$. By definition, an {\em extended standard basis} is a list of $1$-forms
$$
\omega_{-1},\omega_0,\omega_1,\ldots,\omega_{s+1},
$$
such that $\nu_{\mathcal C}(\omega_i)=\lambda_i$ for $i=-1,0,1,\ldots,s$ and $\mathcal C$ is an invariant branch of $\omega_{s+1}$, that is $\nu_{\mathcal C}(\omega_{s+1})=\infty$, with some restrictions on the weighted order of $\omega_{s+1}$.

Associated to the final divisor $E$ given by $\mathcal C$, we have a {\em divisorial order} $\nu_E(\omega)$ defined for functions and $1$-forms. In adapted coordinates it is the weighted monomial order that assigns the weight $an+bm$ to the monomial $x^ay^b$. Both the differential values and the divisorial orders act ``like'' valuations and we have that $\nu_E(\omega)\leq\nu_{\mathcal C}(\omega)$. For the case of a function we have that if $\nu_E(df)<nm$, then there is no resonance in the sense that $\nu_E(df)=\nu_{\mathcal C}(d f)$. Thus, the ``new differential values'' in $\Lambda^{\mathcal C}$ will correspond to resonant $1$-forms $\omega$ such that $\nu_E(\omega)<\nu_{\mathcal C}(\omega)$.

The structure of the semimodule $\Lambda^\mathcal C$ is well known (see \cite{delorme,patricio, patriciomoyano}); anyway, we provide complete proofs using another approach  in the appendices of the paper. The key elements are the {\em axes} $u_i$, and the {\em critical orders $t_i$}, defined by
$$
u_{i+1}=\min (\Lambda_{i-1}\cap (\lambda_i+\Gamma)),\quad t_{i+1}=t_i+u_{i+1}-\lambda_i,
$$
starting at $u_0=n$ and $t_{-1}=n,t_0=m$, where $\Lambda_{i-1}^{\mathcal{C}}=\cup_{j=-1}^{i-1}(\lambda_j+\Gamma)$. The axes are defined for $i=0,1,\ldots,s+1$ and the critical orders for $i=-1,0,\ldots,s+1$. We know that the semimodule is {\em increasing } in the sense that $\lambda_i>u_i$ for $i=1,2,\ldots,s$ and the elements of any extended standard basis are characterized by the following properties
\begin{enumerate}
\item $\nu_E(\omega_i)=t_i$ and
$\nu_\mathcal C(\omega_i)\notin \Lambda_{i-1}^\mathcal C$, for  $i=-1,0,\ldots,s$.
    \item $\nu_E(\omega_{s+1})=t_{s+1}$ and $\nu_\mathcal C(\omega_{s+1})=\infty$.
\end{enumerate}
Of course, the above properties assure that $\nu_{\mathcal C}(\omega_i)=\lambda_i$.

From the geometrical viewpoint, for each $i=1,2,\ldots,s+1$, the elements $\omega_i$  of an extended standard basis are what we call {\em basic and resonant}. This property implies that the transform $\tilde \omega_i$ of the $1$-form $\omega_i$ by the morphism $\pi$ of reduction of singularities of $\mathcal C$ has two remarkable properties:
\begin{enumerate}
\item[a)] The greatest common divisor  of the coefficients of $\tilde\omega_i$ defines a normal crossings divisor at the points of $E$ contained in the exceptional divisor of the morphism $\pi$.
    \item[b)] The divisor $E$ is dicritical (not invariant) for the foliation given by $\tilde\omega_i=0$. Moreover, this foliation is nonsingular and it has normal crossings with the exceptional divisor of $\pi$ at the points of $E$.
\end{enumerate}
As a consequence of this, given an extended standard basis, we find a {\em dicritical package $\{\mathcal C^i_P\}$} of cusps for each $i=1,2,\ldots,s+1$ parameterized by the points $P\in E$ that are not corners of the exceptional divisor (that is, elements of $\mathbb C^*$). Each $\mathcal C^i_P$ corresponds to the invariant curve of $\tilde\omega_i=0$ through the point $P$. In particular, if $P_0$ is the infinitely near point of $\mathcal C$ at $E$, we have that $\mathcal C^{s+1}_{P_0}=\mathcal C$. In a terminology inspired in Equisingularity Theory and Reduction of Singularities (see for instance \cite{Abh-M, wall} and \cite{seidenberg} for the case of foliations), we could say that $\{\mathcal C^i_{P_0}\}$ are the {\em specific analytic semiroots} and that $\{\mathcal C^i_{P}\}$ are the {\em general analytic semiroots} of $\mathcal C$ associated to the given extended standard basis.

The property of $E$ to be dicritical for the $1$-forms $\omega_i$ has been
suggested to us by M. E. Hernandes. We have a work in progress with him in this direction \cite{Cor-H-H}.

The main objective of this paper is to describe the semimodule and extended standard bases of the analytic semiroots. The statement is the following one:
\begin{theorem} Let $\Lambda^\mathcal{C}=\cup_{j=-1}^{s}(\lambda_i+\Gamma)$ be the semimodule of differential values  and consider an extended standard basis
$$
\omega_{-1}=dx,\omega_0=dy,\omega_1,\ldots,\omega_{s+1}
$$
of the cusp $\mathcal C$. Take an index $i\in\{1,2,\ldots,s+1\}$ and an analytic semiroot $\mathcal C^{i}_P$ of $\mathcal C$ associated to the given extended standard basis. Then the semimodule of differential values of $\mathcal C^i_P$ is precisely $\Lambda_{i-1}^\mathcal{C}$ and
$$
\omega_{-1}=dx,\omega_0=dy,\omega_1,\ldots,\omega_{i}
$$
is a extended standard basis for $\mathcal C^i_P$.
\end{theorem}

The proof of this result uses as a main tool Delorme's decomposition of the elements of a standard basis. In the appendices, we provide  proofs, using a different approach  to the one of Delorme, of the structure results for the semimodule of differential values and of Delorme's decomposition.

%

Let us remark that it is possible to have curves of the dicritical package of the elements $\omega_j$, when $j\geq 2$, of an extended standard basis that are not analytically equivalent, although they have the same semimodule of differential values. This occurs for instance if we compute a standard basis for the curve
$$
t\mapsto (t^7, t^{17}+t^{30}+t^{33}+t^{36}).
$$
as shown in Example~\ref{ex:7-17}. A natural question arises about ``how many'' analytic classes may be obtained  in this way.

\medskip
\noindent{\bf Acknowledgements:} We are grateful to Marcelo E. Hernandes for all the conversations and suggestions on the subject. The work of Oziel Gómez has influenced our work in this paper, as well as the guidance of Patricio Almirón on the study of semimodules.
\section{Cusps and Cuspidal Divisors}
We are interested in the analytic moduli of branches with only one Puiseux pair, the {\em analytic cusps}. The last divisor of the minimal reduction of singularities of an analytic cusp is what we call a {\em cuspidal divisor}. As we shall see below, the study of the analytic moduli may be done through a fixed cuspidal divisor.
\subsection{Cuspidal sequences of blowing-ups}
Our ambient space is a two-dimen\-sional germ of nonsingular complex analytic space $(M_0,P_0)$. We are going to consider a specific type of finite sequences of blowing-ups centered at points, that we call {\em cuspidal sequences of blowing-ups} and we introduce below.

First of all, let us establish some notations concerning a nonempty finite sequence of blowing-ups
centered at points
$$
\mathcal S=\{\pi_k:(M_k,K_k)\rightarrow (M_{k-1},K_{k-1});\;  k=1,2,\ldots,N\},
$$
starting at $(M_0,P_0)=(M_0,K_0)$. For any $k=1,2,\ldots,N$, the center of $\pi_k$ is denoted by $P_{k-1}$, note that $P_{k-1}\in K_{k-1}$. We denote the intermediary morphisms as
$\sigma_k:(M_k,K_k)\rightarrow (M_0,P_0)$ and $\rho_k: (M_N,K_N)\rightarrow (M_k,K_k)$, where
$$
\sigma_k=\pi_1\circ\pi_2\circ\cdots\circ\pi_k,\quad
\rho_k=\pi_{k+1}\circ\pi_{k+2}\circ\cdots\circ\pi_N.
$$
We denote the exceptional divisor of $\pi_k$ as $E^k_k=\pi_{k}^{-1}(P_{k-1})$.
 By induction, for any $1\leq j<k$ we denote by $E^k_j\subset M_k$ the strict transform of $E^{k-1}_j$ by $\pi_k$. In this way we have that
$$
K_k =\sigma_{k}^{-1}(P_0)=E^k_1\cup E^k_2\cup\cdots\cup E^k_k.
$$
For any $P\in K_k$, we define $e(P)=\#\{j; P\in E^k_j\}$. Note that $e(P)\in \{1,2\}$. If $e(P)=1$, we say that $P$ is a {\em free point} and if $e(P)=2$ we say that it is a {\em corner point}. Note that all the points in $E^1_1=K_1$ are free points. The last divisor $E^N_N$ will be denoted $E=E^N_N$. We will also denote $M=M_N$, $K=K_N$ and $\pi=\sigma_N:(M,K)\rightarrow (M_0,P_0)$.
\begin{defn}
Following usual Hironaka's terminology, we say that the sequence $\mathcal S$ is a {\em bamboo} if $P_k\in E^k_k$ for any $k=1,2,\ldots,N-1$. We say that $\mathcal S$ is a {\em cuspidal sequence} if it is a bamboo and $e(P_{k-1})\leq e(P_{k})$, for any $2\leq k\leq N-1$. The last divisor $E$ of a cuspidal sequence is called a {\em cuspidal divisor}.
\end{defn}
\begin{remark} In the frame of Algebraic Geometry, the cuspidal divisor $E$ corresponds to a valuation $\nu_E$ of the field of rational functions and it determines completely the cuspidal sequence, once the starting ambient space is fixed. We will work with this valuation, but we present it in a direct way.
\end{remark}

Given a cuspidal sequence $\mathcal S$ with $N\geq 2$, there is well defined {\em index of freeness}  $f$ with $1\leq f \leq N-1$ such that
$
P_1,P_2,\ldots, P_f
$
are free points and $P_{f+1},P_{f+2},\ldots,P_{N-1}$ are corner points. If $N=1$ we put $f=0$. A nonsingular branch $(Y,P_0)\subset (M_0,P_0)$ has {\em maximal contact with $\mathcal S$} if and only if $P_k$ is an infinitely near point of $(Y,P_0)$ for each $k=1,2,\ldots,P_f$.

\begin{remark}For any cuspidal sequence $\mathcal S$ there is at least one nonsingular branch $(Y,P_0)$ having maximal contact with $\mathcal S$. Moreover, if $(Y,P_0)$ has maximal contact with $\mathcal S$ and $(Y',P_0)$ is another nonsingular branch, we have that $(Y',P_0)$ has maximal contact with $\mathcal S$ if and only if
$
i_{P_0}(Y,Y')\geq f+1
$,
where $i_{P_0}(Y,Y')$ stands for the intersection multiplicity.
\end{remark}

We define intermediate cuspidal sequences of a cuspidal sequence $\mathcal S$ as follows.
Given an index $0\leq j\leq N-1$, the {\em intermediate $j^\text{th}$-cuspidal sequence $\mathcal S^{(j)}$ of $\mathcal S$} is the sequence of length $N-j$, starting at $(M_j,P_j)$  such that the blowing ups
$$
\pi_k^{(j)}: (M_{k+j},K_{k}^{(j)})\rightarrow (M_{k+j-1}, K^{(j)}_{k-1}),\quad  k=1,2,\ldots,N-j
$$
are obtained by restriction from $\pi_{k+j}$, where we put $K^{(j)}_0=\{P_j\}$ and $K_{k}^{(j)}\subset K_{k+j}$ is the image inverse of $P_j$ by $\pi_{j+1}\circ\pi_{j+2}\circ \cdots\circ \pi_{j+k}$.

\begin{remark} Note that the $(k,i)$-divisor of $\mathcal S^{(j)}$ corresponds to the $(k+j,i+j)$ divisor of $\mathcal S$. In particular the last divisors of $\mathcal S^{(j)}$ and $\mathcal S$ are both equal to $E$.
\end{remark}

The {\em Puiseux pair  $(n,m)$} of $\mathcal S$ is defined by an inductive process that corresponds to Euclides' algorithm as follows. If $N=1$, we put $(n,m)=(1,1)$. If $N>1$, we consider the intermediate cuspidal sequence $\mathcal S^{(1)}$ starting at $(M_1,P_1)$ that is supposed to have Puiseux pair $(n_1,m_1)$. Then
    \begin{enumerate}
    \item If $f\geq 2$, we have that $f_1=f-1$ and we put $(n,m)=(n_1,m_1+n_1)$.
    \item If $f=1$, we put $(n,m)=(m_1,n_1+m_1)$.
    \end{enumerate}

 We see that $1\leq n\leq m$ and $n,m$ are without common factor. Note also that $f\geq 2$ if and only if $m\geq 2n$. Moreover, if $f=1$ and $N\geq 2$, we have that $2\leq n<m<2n$.

\begin{prop}
 \label{prop:cuspidadsequence}
 Consider  $1\leq n\leq m$ without common factor and a nonsingular branch $(Y,P_0)\subset (M_0,P_0)$. There is a unique cuspidal sequence $\mathcal S$ starting at $(M_0,P_0)$ having maximal contact with $(Y,P_0)$ and such that $(n,m)$ is the Puiseux pair of  $\mathcal S$.
\end{prop}
\begin{proof} If $n=m=1$, the only possibility is that $N=1$ and then $\mathcal S$ consists in the  blowing-up of $P_0$.  Let us proceed by induction on $n+m$ and assume that $n+m>2$. We necessarily have that $N\geq 2$, the first blowing-up $\pi_1$ is centered in $P_0$ and $P_1$ is the infinitely near point of $Y$ in $E^1_1$.

Assume first that $2n\leq m$. We apply induction to $(Y_1,P_1)$ with respect to the pair $n',m'$ where
$n'=n$, $m'=m-n$
and we obtain a cuspidal sequence $\mathcal S'$ over $(M_1,P_1)$ of length $N'$ with the required properties. We construct $\mathcal S$ of length $N=N'+1$ by taking $\pi_k$ centered at the point $P'_{k-2}$, for $k=2,3,\ldots,N'+1$.

In the case that $n\leq m<2n$ we consider the branch $(Y'_1,P_1)=(E^1_1,P_1)$, we apply
induction to $(Y'_1,P_1)$ with respect to the pair $n',m'$ where
$n'=m-n$,   $m'=n$ and we obtain a cuspidal sequence $\mathcal S'$ over $(M_1,P_1)$ of length $N'$. We construct $\mathcal S$ os length $N=N'+1$ as before.

The uniqueness of $\mathcal S$ follows by an inductive invoking of the uniqueness after one blowing-up.
\end{proof}

We denote by $\mathcal S^{n,m}_Y$ the sequence obtained in Proposition \ref{prop:cuspidadsequence}. Recall that $Y'$ has maximal contact with $\mathcal S^{n,m}_Y$ if and only if $
i_{P_0}(Y,Y')\geq f+1,
$ and hence in this case we have that  $\mathcal S^{n,m}_Y=\mathcal S^{n,m}_{Y'}$. Note also that given a cuspidal sequence $\mathcal S$ there is a nonsingular branch $(Y,P_0)$ and a Puiseux pair $(n,m)$ in such a way that $\mathcal S=\mathcal S^{n,m}_Y$.
\subsection{Coordinates Adapted to a Cuspidal Sequence of Blowing-ups} Consider a cuspidal sequence $\mathcal S$ over $(M_0,P_0)$. A system $(x,y)$ of local coordinates at $P_0$ is {\em adapted to $\mathcal S$} if and only if $y=0$ has maximal contact with $\mathcal S$. In particular, we have that $\mathcal S=\mathcal S^{n,m}_{y=0}$, where $(n,m)$ is the Puiseux pair of $\mathcal S$.

The blowing-ups of $\mathcal S$ have a monomial expression in terms of adapted coordinates as we see below. Assume that $\mathcal S=\mathcal S^{n,m}_{y=0}$, with $N\geq 2$. Let us describe a local coordinate system $(x_1,y_1)$ at $P_1$ and a pair $(n_1,m_1)$:
\begin{itemize}
 \item If $f\geq 2$, we know that $2n\leq m$ and we put
$$
n_1=n,\quad m_1=m-n, \quad x=x_1, \quad y=x_1y_1.
$$
\item
If $f=1$, we have that $2n>m>n\geq 2$ and  we put
$$
n_1=m-n,\quad m_1=n, \quad y=x_1y_1, \quad x=y_1.
$$
\end{itemize}
 The reader can verify that $(x_1,y_1)$ is a coordinate system adapted to $\mathcal S^{(1)}$ and that $(n_1,m_1)$ is its  Puiseux pair. In this way, we have local coordinates $x_j,y_j$ at each $P_j$, for $0\leq j\leq N-1$.

Once we have an adapted coordinate system $(x,y)$, we denote $(H_0,P_0)$ the normal crossings germ given by $xy=0$. Define $H_j=\sigma_j^{-1}(H_0)$, then the germ of $H_j$ at $P_j$ is given by $x_jy_j=0$, for any $0\leq j\leq N-1$. We can also consider $H=\pi^{-1}(H_0)\subset M$; it is a normal crossings divisor of $(M,K)$ containing $K$.

\subsection{Cuspidal Analytic Module} Consider a cuspidal sequence $\mathcal S$ with Puiseux pair $(n,m)$ with $2\leq n$. Let $E$ be the last divisor of $\mathcal S$. We say that an analytic branch $(\mathcal C,P_0)\subset (M_0,P_0)$ is an {\em $E$-cusp}, or a {\em $\mathcal S$-cusp} if the strict transform of $(\mathcal C,P_0)$ under the sequence of blowing-ups $\pi$ is nonsingular and cuts transversely $E$ at a free point. Let us denote by
$
\operatorname{Cusps}(E)= \operatorname{Cusps}(\mathcal S)
$
the family of $E$-cusps.

Each element of  $\operatorname{Cusps}(\mathcal S)$ is equisingular to the irreducible cusp $y^n-x^m=0$, where $(n,m)$ is the Puiseux pair of $\mathcal S$.  Moreover, we have the following result
\begin{prop} Consider a cuspidal sequence $\mathcal S$ with Puiseux pair $(n,m)$ and last divisor $E$.  Let
$({\mathcal C},P_0)$ be a branch of $(M_0,P_0)$ equisingular to the irreducible cusp $y^n-x^m=0$. There is an $E$-cusp analytically equivalent to $({\mathcal C},P_0)$.
\end{prop}
\begin{proof} Choose a local coordinate system $(x,y)$ adapted to $\mathcal S$.

If $n=1$, the branch $({\mathcal C},P_0)$ is nonsingular. Then, there is an automorphism $\phi:(M_0,P_0)\rightarrow (M_0,P_0)$ such that $\phi({\mathcal C})=(y=0)$. We are done since in this case $y=0$ is an $E$-cusp.

 Assume that $2\leq n<m$.  In view of the classical arguments of Hironaka (see for instance~\cite{Ar-H-V}, p. 105), there is a nonsingular branch $(Z,P_0)$ having maximal contact with $({\mathcal C},P_0)$, that is with the property that
$$
i_{P_0}(Z,{\mathcal C})= m.
$$
Take an automorphism $\phi:(M_0,P_0)\rightarrow (M_0,P_0)$ such that $\phi(Z)=(y=0)$. We have that $(\phi(\mathcal C),P_0)$ is an $E$-cusp.
\end{proof}
According to the above result, the analytic moduli of the family of branches equisingular to the irreducible cusp $y^n-x^m=0$ is faithfully represented by the analytic moduli of the family $\operatorname{Cusps}(\mathcal S)$.
\begin{quote}
Along the rest of this paper, we consider a fixed cuspidal sequence $\mathcal S$ where $(n,m)$ is its Puiseux pair and $E$ is the last divisor. Recall also that the composition of all the blowing-ups of $\mathcal S$ is denoted by
$$
\pi:(M,K)\rightarrow (M_0,P_0).
$$
We also choose a local coordinate system $(x,y)$ adapted to $\mathcal S$.
\end{quote}
\section{Divisorial Order}
 Consider a holomorphic function
$h$ in $(M,K)$ defined globally in $E\subset K$, the {\em divisorial order $\nu_E(h)$ of $h$} is obtained as follows. Take a point $P\in E$  and choose a reduced local equation $u=0$ of the germ $(E,P)$, then
$$
\nu_E(h)=\max\{a\in \mathbb Z;\; u^{-a}h\in \mathcal O_{M,P}\}.
$$
This definition does not depend on the chosen point $P\in E$ nor on the local reduced equation of $E$.
 Take a point $P_j$, with $j\in \{0,1,\ldots,N-1\}$
and a germ of holomorphic function $h\in \mathcal O_{M_j,P_j}$. Then $\rho_j^*h$ is a germ of function in $(M,K)$ globally defined in $E$. We define the {\em divisorial order $\nu_E(h)$} by
$
\nu_E(h)=\nu_E(\rho_j^*h)
$.
\begin{prop}\label{prop:valoracionmonomial} Consider a germ $h\in \mathcal O_{M_0,P_0}$ that we write as
$$
h=\sum_{\alpha,\beta}h_{\alpha,\beta}x^\alpha y^\beta,\quad h_{\alpha\beta}\in \mathbb C.
$$
Then $\nu_E(h)=\min\{n\alpha+m\beta;\;  h_{\alpha,\beta}\ne 0\}$.
\end{prop}
\begin{proof} If $n=m=1$ we have a single blowing-up and we recover the usual multiplicity, that we visualize in $E$ as $\nu_E(h)$. Let us work by induction on $n+m$ and assume that $n+m\geq 2$. We remark that
$
\nu_E(h)=\nu_E(\pi_1^*h)
$.
Consider the first intermediate sequence $\mathcal S^{(1)}$ of $\mathcal S$, with adapted coordinates $(x_1,y_1)$. Recalling how we obtain intermediate coordinate systems, we conclude that
\begin{eqnarray*}
\pi_1^*h= \sum_{\alpha,\beta}h_{\alpha,\beta}x_1^{\alpha+\beta} y_1^\beta; \text{ if } f\geq 2; \text{ here } n_1=n, m_1=m-n,\\
\pi_1^*h= \sum_{\alpha,\beta}h_{\alpha,\beta}x_1^{\beta} y_1^{\alpha+\beta}; \text{ if } f=1; \text{ here } n_1=n-m, m_1=n.
\end{eqnarray*}
We end by applying induction hypothesis.
\end{proof}

\subsection{Divisorial Order of a Differential Form}
Recall that we denote
$$H_0=(xy=0)\subset M_0, \quad H_j=\sigma_j^{-1}(H_0)\subset M_j$$ and that $H_j$ is locally  given at $P_j$ by $x_iy_j=0$ for $0\leq j\leq N-1$. We also consider $H_N=H=\pi^{-1}(H_0)\subset M$. Each $H_j$ is a normal crossings divisor in $(M_j,K_j)$, containing $K_j$.

Take a point $Q\in K_j$, not necessarily equal to $P_j$, in particular we consider also the case $j=N$. Select local coordinates $(u,v)$ such that $(u=0)\subset H_j\subset (uv=0)$, then we have that either $H_j=(u=0)$ or $H_j=(uv=0)$ locally at $Q$. The $\mathcal O_{M_j,Q}$-module
$
\Omega^1_{M_j,Q}[\log H_j]
$
of germs of {\em $H_j$-logarithmic 1-forms} is the rank two free  $\mathcal O_{M_j,Q}$-module generated by
\begin{eqnarray*}
du/u, dv&\mbox{ if }& H_j=(u=0),\\
du/u, dv/v&\mbox{ if }& H_j=(uv=0).
\end{eqnarray*}
Note that  $\Omega^1_{M_j,Q}\subset \Omega^1_{M_j,Q}[\log H_j]$. Indeed, a differential $1$-form $\omega=adu+bdv$ may be written as
$$
\omega=ua\frac{du}{u}+bdv=ua\frac{du}{u}+vb\frac{dv}{v}.
$$

Now, let us consider a $1$-form $\omega\in \Omega^1_{M}[\log H]$ defined in the whole divisor $E$ (we suppose that the reader recognizes the sheaf nature of $\Omega^1_{M}[\log H]$). Select a point $Q\in E$ and a local reduced equation $u=0$ of $E$ at $Q$. We define the {\em divisorial order $\nu_E(\omega)$} by
$$
\nu_E(\omega)=\max\{\ell\in \mathbb Z;\; u^{-\ell}\omega\in \Omega^1_{M,Q}[\log H]\}.
$$
The definition is independent of $Q\in E$ and of the reduced local equation of $E$.
\begin{remark} Let $\omega\in \Omega^1_{M}[\log E]$ be globally defined on $E$ as before. Since $E$ is one of the irreducible components of $H$, we have that
$$\Omega^1_{M}[\log E]\subset \Omega^1_{M}[\log H].
$$
Let us choose a reduced local equation $u=0$ of $E$ at a point $Q\in E$ as before. A direct computation shows that
\begin{equation}
\nu_E(\omega)= \max\{\ell\in \mathbb Z;\; u^{-\ell}\omega\in \Omega^1_{M,Q}[\log E]\}.
\end{equation}
This remark shows that the divisorial order, applied to $1$-forms $\omega\in \Omega^1_{M}[\log E]$ is independent of the choice of the adapted coordinate system that defines $H_0$. Anyway, this is only a remark for the case $n=1$, since when $n\geq 2$ the divisor $H$ at the points of $E$ is itself independent of the adapted coordinate system.
\end{remark}

\begin{defn} For any $\omega\in \Omega^1_{M_j,P_j}$, the {\em divisorial order $\nu_E(\omega)$} is defined by
$
\nu_E(\omega)=\nu_E(\rho_{j}^*\omega).
$
\end{defn}

\begin{prop}Consider a differential $1$-form
$
\omega=adx+bdy\in \Omega^1_{M_0,P_0}
$,
that we can write as
$$
\omega=xa({dx}/{x})+yb({dy}/{y})\in \Omega^1_{M_0,P_0}[\log H_0].
$$
Then, we have that
$\nu_E(\omega)=\min\{\nu_E(xa),\nu_E(yb)\}
$.
\end{prop}
\begin{proof} Write $\omega=f(dx/x)+g(dy/y)$. We proceed by induction on $N$.
If $N=1$ we have that
$E=(x'=0)$ where $x=x',y=x'y'$ and
$$
\pi_1^*\omega=(f+g)(dx'/x')+g(dy'/y').
$$
Then
$
\nu_E(\omega)= \min\{\nu_E(f+g),\nu_E(g)\}=
\min\{\nu_E(f),\nu_E(g)\}
$
and we are done. If $N\geq 2$, we have that
$$
\nu_E(\omega)=\nu_E(\pi^*\omega)= \nu_E(\rho_1^*(\pi_1^*\omega))= \nu_E(\pi_1^*\omega).
$$
By induction hypothesis, we have $$
\nu_E(\pi_1^*\omega)= \min\{\nu_E(f+g),\nu_E(g)\}=
\min\{\nu_E(f),\nu_E(g)\}
$$ and we are done as before.
\end{proof}

\begin{corolario} If $f\in \mathcal O_{M_0,P_0}$ and $\omega=df$, then
$
\nu_E(\omega)=\nu_E(f)
$.
\end{corolario}
\begin{proof}
It is enough to write
$
df=x({\partial f}/{\partial x})({dx}/{x})+
y({\partial f}/{\partial y})({dy}/{y}),
$
recalling Euler's identity $gP=xP_x+yP_y$ for degree $g$ homogeneous polynomials.
\end{proof}
\subsection{Weighted Initial Parts}\label{sec:weighted-initial-part}.
Consider a nonzero germ $h\in \mathcal O_{M_0,P_0}$, that we write as
$
h=\sum_{\alpha,\beta}h_{\alpha\beta}x^\alpha y^\beta
$.
Suppose that $q\leq \nu_E(h)$. We define the {\em weighted initial part $\operatorname{In}^q_{n,m; x,y}(h)$} by
$$
\operatorname{In}^q_{n,m; x,y}(h)=\sum_{n\alpha+m\beta=q}h_{\alpha\beta}x^\alpha y^\beta.
$$
Note that $\operatorname{In}^q_{n,m; x,y}(h)=0$ if and only if $q<\nu_E(h)$. Anyway, we can write
$$
h=\operatorname{In}^q_{n,m; x,y}(h)+\tilde h,\quad \nu_E(\tilde h)>q.
$$
This definition extends to logarithmic differential $1$-forms $\omega\in \Omega^1_{M_0,P_0}[\log (xy=0)]$ as follows. Take $q\leq \nu_E(\omega)$. Write
$
\omega=f(dx/x)+g(dy/y)
$. We define
$$
\operatorname{In}^q_{n,m; x,y}(\omega)= \operatorname{In}^q_{n,m; x,y}(f)(dx/x)+\operatorname{In}^q_{n,m; x,y}(g)(dy/y).
$$
As before, we have
$
\omega= \operatorname{In}^q_{n,m; x,y}(\omega)+\tilde \omega$, with  $\nu_E(\tilde \omega)>q
$.
\begin{remark} The definition of initial part we have presented should be made in terms of graduated rings and modules to be free of coordinates. Anyway, this ``coordinate-based'' definition is enough for our purposes.
\end{remark}
\begin{prop} Assume that $N>1$, take $\omega\in \Omega^1_{M_0,P_0}[\log (xy=0)]$ and $q\in \mathbb Z_{\geq 0}$ with $q\leq\nu_E(\omega)$. If
$
W=\operatorname{In}^q_{n,m; x,y}(\omega)
$, then $\pi_1^*(W)=\operatorname{In}^q_{n_1,m_1; x_1,y_1}(\pi_1^*\omega)$.
\end{prop}
\begin{proof} Left to the reader.
\end{proof}
\section{Total Cuspidal Dicriticalness}
This section is devoted to characterize the $1$-forms $\omega\in \Omega^1_{M_0,P_0}$ whose transform $\pi^*\omega$ defines a foliation that is transversal to $E$ and has normal crossings with $K$ at  each point of $E$. These 1-forms are the so-called pre-basic and resonant 1-forms. We detect these properties in terms of resonances of the initial part. The initial part is visible in the Newton polygon as the contribution of the 1-form to a single vertex $(a,b)$, under the condition that  the Newton polygon is contained in the particular region $R^{n,m}(a,b)$.

\subsection{Reduced Divisorial Order and Basic Forms} Let us consider a nonnull differential $1$-form $\omega\in \Omega^1_{M_0,P_0}$. Let $V_\omega=x^a y^b$ be the monomial defined by the property that
$
\omega=V_\omega\eta
$,
where $\eta\in \Omega^1_{M_0,P_0}[\log (xy=0)]$ is a logarithmic form that cannot be divided by any nonconstant monomial.
 We define the {\em reduced divisorial order $\operatorname{rdo}_E(\omega)$} to be
$
\operatorname{rdo}_E(\omega)=\nu_E(\eta)
$.
\begin{defn} We say that $\omega\in \Omega^1_{M_0,P_0}$ is a {\em basic $1$-form} if and only if its reduced divisorial order satisfies that $\operatorname{rdo}_E(\omega)<nm$.
\end{defn}

\begin{prop}
 \label{prop:stabilityofbasic}
 Assume that $N\geq 2$ and take $\omega\in \Omega^1_{M_0,P_0}$. If $\omega$ is a basic $1$-form, then $\pi_1^*\omega$ is also a basic $1$-form.
\end{prop}
\begin{proof} Put $p=\operatorname{rdo}_E(\omega)=\nu_E(\eta)<nm$. Recall that $\nu_E(\eta)=\nu_E(\pi_1^*\eta)$. Since monomials are well behaved under $\pi_1$, it is enough to show that there are $c,d\geq 0$ such that
$
\pi_1^*\eta=x_1^{c}y_1^{d}\eta'$, with
$\nu_E(\eta')<n_1m_1 $.
 Write
$$
\eta=\sum_{\alpha,\beta} x^{\alpha}y^{\beta}\eta_{\alpha\beta},\quad  \eta_{\alpha\beta}=
\mu_{\alpha\beta}\frac{dx}{x}+\zeta_{\alpha\beta}\frac{dy}{y}, \quad (\mu_{\alpha\beta},\zeta_{\alpha\beta})\in \mathbb C^2.
$$
Recall that
$
p=\min\{n\alpha+m\beta;\; \eta_{\alpha\beta}\ne 0\}
$.
Put $r=\min\{\alpha+\beta;\; \eta_{\alpha\beta}\ne 0\}$.
We have two cases: $f=1$ and $f\geq 2$, where $f$ is the index of freeness.

Assume first that $f\geq 2$ and hence $2n\leq m$. In this situation, we have that $x=x_1$, $y=x_1y_1$, $n_1=n$, $m_1=m-n\geq n$ and
$$
\pi_1^*(\eta)=x_1^r\eta',\quad \eta'=\sum_{\alpha,\beta}x_1^{\alpha+\beta-r}y_1^\beta\eta'_{\alpha\beta},\quad
\eta'_{\alpha\beta}=\left(\mu_{\alpha\beta}+\zeta_{\alpha\beta}\right)\frac{dx_1}{x_1}
+\zeta_{\alpha\beta}\frac{dy_1}{y_1}.
$$
Note that $\eta'_{\alpha\beta}\ne 0$ if and only if $\eta_{\alpha\beta}\ne 0$. Hence
\begin{eqnarray*}
\nu_E(\eta')&=&\min \{n_1(\alpha+\beta-r)+m_1\beta;\eta_{\alpha\beta}\ne0\}=\\
&=& \min \{n
(\alpha+\beta-r)+(m-n)\beta;\eta_{\alpha\beta}\ne 0
\}=\\
&=& \min \{n\alpha+m\beta-nr;\eta_{\alpha\beta}\ne 0\}=p-nr.
\end{eqnarray*}
We have to verify that $p-nr<n_1m_1$, where $n_1m_1=n(m-n)=nm-n^2$. If $r\geq n$, we are done, since by hypothesis we have that $p<nm$. Assume that $r<n$. There are $\tilde\alpha,\tilde\beta$ with $\eta_{\tilde\alpha\tilde\beta}\ne 0$ such that $\tilde\alpha+\tilde\beta=r$. Then
\begin{eqnarray*}
p-nr&\leq&
n\tilde\alpha+m\tilde\beta-nr=n(\tilde\alpha+\tilde\beta)+(m-n)\tilde\beta-nr=\\
&=&(m-n)\tilde\beta<(m-n)n,
\end{eqnarray*}
since $\tilde\beta\leq r<n$.

Assume that $f=1$ and thus $n<m<2n$. We have $x=y_1$, $y=x_1y_1$, $n_1=m-n<n$, $m_1=n$ and
$$
\pi_1^*(\eta)=y_1^r\eta'',\quad \eta''=\sum_{\alpha,\beta}x_1^{\beta}y_1^{\alpha+\beta-r}\eta''_{\alpha\beta},\quad
\eta''_{\alpha\beta}=\zeta_{\alpha\beta}\frac{dx_1}{x_1}
+\left(\mu_{\alpha\beta}+\zeta_{\alpha\beta}\right)\frac{dy_1}{y_1}.
$$
As before, we have that $\eta''_{\alpha\beta}\ne 0$ if and only if $\eta_{\alpha\beta}\ne 0$. Hence
\begin{eqnarray*}
\nu_E(\eta'')&=&\min \{n_1\beta+m_1(\alpha+\beta-r);\; \eta_{\alpha\beta}\ne0\}=\\
&=&
\min \{(m-n) \beta+n(\alpha+\beta-r);\; \eta_{\alpha\beta}\ne 0\}=
\\
&=& \min \{m\beta+n\alpha-nr;\; \eta_{\alpha\beta}\ne 0\}=p-nr.
\end{eqnarray*}
We verify that $p-nr<n_1m_1$ exactly as before.
\end{proof}

\subsection{Resonant Basic forms} Let $\omega\in \Omega^1_{M_0,P_0}$ be a basic 1-form with $p=\operatorname{rdo}_E(\omega)$. This means that there is $\eta\in \Omega^1_{M_0,P_0}[\log (xy=0)]$ and $a,b\geq 0$ such that $\omega=x^ay^b\eta$
$$
\omega=x^ay^b\eta, \quad \eta\in \Omega^1_{M_0,P_0}[\log (xy=0)],
$$
where $p=\nu_E(\eta)<nm$. The initial part of $\omega$ may be written
$$
\operatorname{In}^{p+na+mb}_{n,m;x,y}(\omega)=x^ay^b W, \quad W=\operatorname{In}^{p}_{n,m;x,y}(\eta).
$$
Note that there is exactly one pair $(c,d)\in \mathbb Z^2_{\geq 0}$ such that $cn+dm=p$.  Then we have that
$$
W=x^cy^d\left\{\mu\frac{dx}{x}+\zeta\frac{dy}{y}\right\}.
$$
We say that $\omega$ is {\em resonant} if and only if $n\mu+m\zeta=0$.

We have the following result that follows directly from the computations in the proof of Proposition \ref{prop:stabilityofbasic}:
\begin{corolario}
 \label{cor:stabilityofresonant}
 Assume that $N\geq 2$.  A basic differential $1$-form $\omega\in \Omega^1_{M_0,P_0}$ is resonant if and only if $\pi_1^*\omega$ is resonant.
\end{corolario}
\subsection{Pre-Basic Forms} Let us introduce a slightly more general class of $1$-forms that we call {\em pre-basic forms}. Given a  $1$-form
\begin{equation}
\label{eq:omegamonomios}
\omega=\sum_{\alpha,\beta}c_{\alpha\beta}x^\alpha y^\beta \omega_{\alpha\beta},\quad \omega_{\alpha\beta}= \left\{\mu_{\alpha\beta}\frac{dx}{x}+\zeta_{\alpha\beta}\frac{dy}{y}\right\},
\end{equation}
the {\em cloud of points $\operatorname{Cl}(\omega;x,y)$} is $\operatorname{Cl}(\omega;x,y)= \{(\alpha,\beta);\omega_{\alpha\beta}\ne 0\}$
and the {\em Newton Polygon $\mathcal N(\omega;x,y)$} is the positive convex hull of
$\operatorname{Cl}(\omega;x,y)$ in $\mathbb R^2_{\geq 0}$.

 Consider a pair $(n,m)$ with $1\leq n\leq m$ such that $n,m$ have no common factor. There are unique $b,d\in \mathbb Z_{\geq 0}$ such that $dn-bm=1$ with the property that $0\leq b<n$ and $0<d\leq m$. We call $(b,d)$ the {\em co-pair} of $(n,m)$.
 \begin{remark}
  \label{rk:thecopair}
  Suppose that $1\leq n\leq  m$ are without common factor and take $b,d$ such that $dn-bm=1$. If $0\leq b<n$, we have that $0<d\leq m$ and then $(b,d)$ is the co-pair of $(n,m)$. In the same way, if $0<d\leq m$, we have that $0\leq b<n$ and then $(b,d)$ is the co-pair of $(n,m)$.
 \end{remark}
\begin{defn}
Given a pair $1\leq n\leq m$ without common factor, we define the region $R^{n,m}$ by  $R^{n,m}= H^{n,m}_-\cap H^{n,m}_+$, where
\begin{eqnarray*}
H^{n,m}_-&=&\{(\alpha,\beta)\in \mathbb R^2;\; (n-b)\alpha+(m-d)\beta\geq 0\},\\
H^{n,m}_+&=&\{(\alpha,\beta)\in \mathbb R^2; \; b\alpha+d\beta\geq 0\},
\end{eqnarray*}
and $(b,d)$ is the co-pair of $(n,m)$.
\end{defn}
\begin{remark} If $n=m=1$, the co-pair of $(1,1)$ is $(b,d)=(0,1)$. Then
$$
H^{1,1}_-=\{(\alpha,\beta);\; \alpha\geq 0\},\quad H^{1,1}_+=\{(\alpha,\beta);\; \beta \geq 0\}.$$
Thus, we have that $R^{1,1}$ is the quadrant $R^{1,1}=\mathbb R^2_{\geq 0}$.
\end{remark}

\begin{remark} The slopes $-(n-b)/(m-d)$ and $-b/d$ satisfy that
$$
-(n-b)/(m-d)<-n/m<-b/d.
$$
Indeed,
we have
$
-n/m<-b/d \Leftrightarrow -dn<-mb=-dn+1
$. On the other hand
\begin{eqnarray*}
-(n-b)/(m-d)<-n/m&\Leftrightarrow&
m(n-b)>n(m-d) \Leftrightarrow bm<dn=bm+1.
\end{eqnarray*}
\end{remark}
We conclude  that
$R^{n,m}$ is a positively convex region of $\mathbb R^{2}$ such that $(0,0)$ is its only vertex  and we have that
$$
R^{n,m}\cap \{(\alpha,\beta)\in \mathbb R^2;\; n\alpha+m\beta=0\}=\{(0,0)\}.
$$
Given a point $(a,b)\in \mathbb R_{\geq 0}^2$, we define $R^{n,m}(a,b)$ by
$
R^{n,m}(a,b)= R^{n,m}+(a,b)
$.

\begin{defn}
We say that $\omega\in \Omega^1_{M_0,P_0}$ is {\em pre-basic} if and only if there is a point $(a,b)\in \operatorname{Cl}(\omega;x,y)$ such that
$\operatorname{Cl}(\omega;x,y)\subset R^{n,m}(a,b)$.
\end{defn}
\begin{remark} Note that $\omega$ is pre-basic if and only if $(a,b)\in \mathcal N(\omega;x,y)$ and $\mathcal N(\omega;x,y)\subset R^{n,m}(a,b)$.
\end{remark}
If $\omega$ is pre-basic,
we have that
$$
\operatorname{Cl}(\omega;x,y) \cap \{(\alpha,\beta)\in \mathbb R^2;\; n\alpha+m\beta=\nu_E(\omega)\}=\{(a,b)\}.
$$
Thus, the initial part $W$ of $\omega$ has the form
\begin{equation}
\label{eq:initialform}
W=x^ay^b\left\{\mu_{ab}\frac{dx}{x}+\zeta_{ab}\frac{dy}{y}\right\}.
\end{equation}
As for basic forms, we say that $\omega$ is {\em resonant} if and only if $n\mu_{ab}+m\zeta_{ab}=0$.
\begin{lema}
\label{lema:prebasicinduccion}
 Assume that $1\leq n<m$, where $n,m$ are without common factor and let $(b,d)$ be the co-pair of $(n,m)$. Let us put
 $(n_1,m_1)=(n,m-n)$, if  $m\geq 2n$ and $(n_1,m_1)=(m-n,n)$, if  $m<2n$.
Then, the co-pair $(b_1,d_1)$ of $(n_1,m_1)$ is given by
$(b_1,d_1)= (b,d-b)$, if $m\geq 2$, and by $(b_1,d_1)=(m-n-d+b,n-b)$, if $m<2n$.
Moreover, we have that $\Psi(R^{n,m})=R^{n_1,m_1}$, where $\Psi$ is
the linear automorphism of $\mathbb R^2$ given by
$\Psi(\alpha,\beta)= (\alpha+\beta,\beta)$, if  $m\geq 2n$, and
$\Psi(\alpha,\beta)= (\beta,\alpha+\beta)$, if $m<2n$.
\end{lema}
\begin{proof} Let us show the first statement.
If $m\geq 2n$,  we have that
$$d_1n_1-b_1m_1= (d-b)n-b(m-n)=1.$$
Moreover, since $0\leq b_1=b<n_1=n$ we conclude that $(b_1,d_1)$ is the co-pair of $(n_1,m_1)$, in view of Remark \ref{rk:thecopair}. If $m<2n$, we have
$$
d_1n_1-b_1m_1= (n-b)(m-n)-(m-n-d+b)n=1.
$$
We know that $0\leq b<n$, hence $0<d_1=n-b\leq m_1=n$; by Remark \ref{rk:thecopair}, we deduce that $(b_1,d_1)$ is the co-pair of $(n_1,m_1)$.

Let us show the second statement.  Consider $(\alpha,\beta)\in \mathbb R^2$ and put
$
(\alpha_1,\beta_1)=\Psi(\alpha,\beta)
$.

{\em Case $m\geq 2n$.}
In order to prove that $\Psi(R^{n,m})=R^{n_1,m_1}$ it is enough to see that
$$
(\alpha,\beta)\in H^{n,m}_-\Leftrightarrow (\alpha_1,\beta_1)\in H^{n_1,m_1}_-\text{ and }
(\alpha,\beta)\in H^{n,m}_+\Leftrightarrow (\alpha_1,\beta_1)\in H^{n_1,m_1}_+.
$$
We verify these properties as follows:
\begin{eqnarray*}
(\alpha_1,\beta_1)\in H^{n_1,m_1}_-&\Leftrightarrow & (n_1-b_1)\alpha_1+(m_1-d_1)\beta_1\geq 0 \Leftrightarrow \\
&\Leftrightarrow & (n-b)(\alpha+\beta)+(m-n-d+b)\beta\geq 0 \Leftrightarrow \\
&\Leftrightarrow & (n-b)\alpha+(m-d)\beta\geq 0 \Leftrightarrow (\alpha,\beta)\in H^{n,m}_-.
\\
(\alpha_1,\beta_1)\in H^{n_1,m_1}_+&\Leftrightarrow & b_1\alpha_1+d_1\beta_1\geq 0 \Leftrightarrow \\
&\Leftrightarrow & b(\alpha+\beta)+(d-b)\beta\geq 0 \Leftrightarrow \\
&\Leftrightarrow & b\alpha+d\beta\geq 0 \Leftrightarrow  (\alpha,\beta)\in H^{n,m}_+.
\end{eqnarray*}

 {\em Case $m<2n$}. In this case, we have that
\begin{eqnarray}
\label{eq:regionR1}
(\alpha,\beta)\in H^{n,m}_+&\Leftrightarrow& (\alpha_1,\beta_1)\in H^{n_1,m_1}_- \\
\label{eq:regionR2}
(\alpha,\beta)\in H^{n,m}_-&\Leftrightarrow& (\alpha_1,\beta_1)\in H^{n_1,m_1}_+.
\end{eqnarray}
and this also implies that $\Psi(R^{n,m})=R^{n_1,m_1}$.
We verify the properties in Equations \eqref{eq:regionR1} and \eqref{eq:regionR2} as follows:
\begin{eqnarray*}
(\alpha_1,\beta_1)\in H^{n_1,m_1}_-&\Leftrightarrow & (n_1-b_1)\alpha_1+(m_1-d_1)\beta_1\geq 0 \Leftrightarrow \\
&\Leftrightarrow & (m-n-(m-n-d+b))\beta+(n-n+b)(\alpha+\beta)\geq 0 \Leftrightarrow \\
&\Leftrightarrow & d\beta+b\alpha\geq 0 \Leftrightarrow (\alpha,\beta)\in H^{n,m}_+.\\
(\alpha_1,\beta_1)\in H^{n_1,m_1}_+&\Leftrightarrow & b_1\alpha_1+d_1\beta_1\geq 0 \Leftrightarrow \\
&\Leftrightarrow & (m-n-d+b)\beta+(n-b)(\alpha+\beta)\geq 0 \Leftrightarrow \\
&\Leftrightarrow & (m-d)\beta+(n-b)\alpha\geq 0 \Leftrightarrow (\alpha,\beta)\in H^{n,m}_-.
\end{eqnarray*}
The proof is ended.
\end{proof}

\begin{prop}
 \label{prop:stabilityoprebasic}
 Assume that $N\geq 2$. For any $\omega\in \Omega^1_{M_0,P_0}$, we have
 \begin{enumerate}
 \item  $\omega$ is pre-basic if and only if $\pi_1^*\omega$ is pre-basic.
 \item $\omega$ is pre-basic and resonant if and only if $\pi_1^*\omega$ is pre-basic and resonant.
 \end{enumerate}
\end{prop}
\begin{proof}We consider two cases as in the statement of Lemma \ref{lema:prebasicinduccion}, the case $m\geq 2n$ and $m<2n$ and we define the linear automorphism $\Psi$ accordingly to these cases, as well as the Puiseux pair $(n_1,m_1)$.  A monomial by monomial computation shows that
\begin{equation}
\label{eq:nubedepuntos}
\operatorname{Cl}(\pi_1^*\omega;x_1,y_1)=\Psi(\operatorname{Cl}(\omega;x,y)).
\end{equation}
In view of Lemma \ref{lema:prebasicinduccion}, we have that
\begin{equation}
\label{eq:regionRab}
\Psi(R^{n,m}(a,b))=R^{n_1,m_1}(\Psi(a,b)).
\end{equation}
Statement (1) is now a direct consequence of Equations \eqref{eq:nubedepuntos} and
\eqref{eq:regionRab}.
Property (2) is left to the reader.
\end{proof}

\begin{prop}
 \label{prop:basicprebasic}
 Take a differential $1$-form $\omega\in \Omega^1_{M_0,P_0}$. We have
 \begin{enumerate}
 \item If $N=1$, then $\omega$ is pre-basic if and only if it is basic.
 \item  If $\omega$ is basic then it is pre-basic.
 \item If $\omega$ is basic and resonant then it is pre-basic and resonant.
 \end{enumerate}
 \end{prop}
 \begin{proof}
 If $N=1$, we have $n=m=1$ and $R^{1,1}(a,b)=\mathbb R^2_{\geq 0}+(a,b)$. Then being basic is the same property of being pre-basic: the Newton Polygon has a single vertex.

 Assume now that $\omega$ is basic. In view of the stability result in Proposition  \ref{prop:stabilityofbasic}, we have that $\tilde\omega$ is basic, where $\tilde\omega$ is the pull-back of $\omega$ in the last center $P_{N-1}$ of the cuspidal sequence. By the previous argument we have that $\tilde\omega$ is pre-basic. Now we apply Proposition  \ref{prop:stabilityoprebasic} to conclude that $\omega$ is pre-basic.

 The resonance for pre-basic $1$-forms that are basic ones is the same property as for basic $1$-forms.
 \end{proof}
\subsection{Totally $E$-dicritical forms} Consider a $1$-form $\omega\in \Omega^1_M$ defined around the divisor $E$. Recall that we have a normal crossings divisor $H$ such that $H\supset E$, coming from our choice of adapted coordinates, although if $n\geq 2$ the divisor $H$ around $E$ is intrinsically defined and it coincides with $K$. We say that $\omega$ is {\em totally $E$-dicritical with respect to $H$} if for any point $P\in E$ there are local coordinates $u,v$ such that $E=(u=0)$, $H\subset (uv=0)$ and $\omega$ has the form
$$
\omega=u^av^bdv,
$$
where $b=0$ when $H=(u=0)$.  Note that $\omega$ defines a non-singular foliation around $E$, this foliation has normal crossings with $H$ and $E$ is transversal to the leaves.

\begin{prop} For any $\omega\in \Omega^1_{M_0,P_0}$, the following properties are equivalent:
\begin{enumerate}
\item $\pi^*\omega$ is totally $E$-dicritical with respect to $H$.
\item The $1$-form $\omega$ is pre-basic and resonant.
\end{enumerate}
\end{prop}
\begin{proof} In view of the stability of the property ``pre-basic and resonant'' under the blowing-ups of $\mathcal S$ given in Proposition~\ref{prop:stabilityoprebasic}, it is enough to consider the case when $N=1$. In this case we have a single blowing-up, $H_0=(xy=0)$ and the property for $\pi^*\omega$ of being totally $E$-dicritical with respect to $H$ is equivalent to say that
$$
\omega=h(x,y)x^a y^b\left[\left\{\frac{dx}{x}-\frac{dy}{y}\right\}+
\sum_{\alpha+\beta\geq 1}x^\alpha y^\beta \left\{\mu_{\alpha\beta}\frac{dx}{x}+\zeta_{\alpha\beta}\frac{dy}{y}\right\}
\right],\quad a,b\geq 1,
$$
where $h(0,0)\ne 0$.
That is, the $1$-form $\omega$ is pre-basic and resonant.
\end{proof}

\begin{remark} If $n\geq 2$ the axes $x'y'=0$ around $P_{N-1}$ coincide with the germ of $K_{N-1}=\sigma_{N-1}^{-1}(P_0)$ at $P_{N-1}$. In this situation, the property of being basic and resonant does not depend on the chosen adapted coordinate system.
\end{remark}

\begin{defn}
 \label{def:omegacusps}
 Given a resonant pre-basic $1$-form $\omega$, we say that a branch $(\mathcal C,0)$ in $(M_0,P_0)$ is a {\em $\omega$-cusp} if and only if it is invariant by $\omega$ and the strict transform of $(\mathcal C,P_0)$ by $\pi$ cuts $E$ at a free point.
\end{defn}
Let us note that each free point of $E$ defines a $\omega$-cusp and conversely, in view of the fact that $\pi^*\omega$ is totally $E$-dicritical with respect to $H$.

One of the results in this paper is that any element of $\operatorname{Cusps}(\mathcal S)$
is a $\omega$-cusp for certain resonant basic $\omega$ and hence can be included in the corresponding ``dicritical package''.

\section{ Differential Values of a Cusp}
 Let us consider a branch $(\mathcal C,P_0)\subset (M_0,P_0)$ belonging to $\operatorname{Cusps}(E)$. It has a Puiseux expansion of the form
$$
(x,y)=\phi(t)=(t^n,\alpha t^m+t^{m+1}\xi(t)),\; \alpha\ne 0.
$$
defined by the fact that for any germ $h\in \mathcal O_{M_0,P_0}$ we have that $(\mathcal C,P_0)\subset (h=0)$ if and only if $h\circ\phi=0$. We recall that the intersection multiplicity of $(\mathcal C, P_0)$ with a germ $h$ is given by
$$
\operatorname{i}_{P_0}(\mathcal C, h)=\operatorname{order}_t(h\circ \phi).
$$
We also denote $\nu_\mathcal C(h)=\operatorname{i}_{P_0}(\mathcal C, h)$.
The {\em semigroup $\Gamma$} of $\mathcal C$ is defined by
$$
\Gamma\cup \{\infty\}=\{\nu_\mathcal C(h);\; h\in \mathcal O_{M_0,P_0}\}.
$$
As stated in Zariski's Equisingularity Theory, this semigroup depends only on the equisingularity class (or topological class) of $\mathcal C$. In our case, we know that all the elements in $\operatorname{Cusps}(E)$ are equisingular to the cusp $y^n-x^m=0$. Hence $\Gamma$ does not depend on the particular choice of $\mathcal C\in \operatorname{Cusps}(E)$. More precisely, we know that $\Gamma$ is the subsemigroup of $\mathbb Z_{\geq 0}$ generated by $n,m$. That is
$$
\Gamma=\{an+bm;\; a,b\in \mathbb Z_{\geq 0}\}.
$$
An important feature of $\Gamma$ is the existence of its conductor $c_\Gamma=(n-1)(m-1)$, which is the smallest element $c_\Gamma \in \Gamma$ such that any non-negative integer greater or equal to $c_\Gamma$ is contained in $\Gamma$.  In a more algebraic way, the conductor ideal $(t^{c_\Gamma})$ is contained in the image of the morphism
$$
\phi^\#: \mathbb C\{x,y\}\rightarrow \mathbb C\{t\},\quad f\mapsto f\circ \phi.
$$

On the other hand, as it was pointed by Zariski, {\em the differential values of $\mathcal C$} may strongly depend on the analytic class of $\mathcal C$. In fact, they are the main discrete invariants in the analytic classification of branches (see \cite{hefez2}).

Given a differential $1$-form $\omega\in \Omega^1_{M_0,P_0}$ with $\omega=gdx+hdy$. If we write $\phi(t)=(x(t),y(t))$, we have that $\phi^*(\omega)=(g(\phi(t))x'(t)+h(\phi(t)) y'(t))dt$. We put $a(t)=t (g(\phi(t))x'(t)+h(\phi(t)) y'(t))$, hence
$$
\phi^*(\omega)=a(t)\frac{dt}{t}
$$
and we define the {\em differential value} $\nu_{\mathcal C}(\omega)$ by
$
\nu_{\mathcal C}(\omega)=\operatorname{order}_t(a(t))
$.

We know that $(\mathcal C,P_0)$ is an {\em invariant branch of $\omega$} if and only if $\phi^*(\omega)=0$ and hence $\nu_{\mathcal C}(\omega)=\infty$. The {\em semimodule $\Lambda^{\mathcal C}$ of the differential values} is defined by
$$
\Lambda^\mathcal C = \{\nu_{\mathcal C}(\omega);\quad \omega \in\Omega^1_{M_0,P_0}, \nu_{\mathcal C}(\omega)\ne \infty \}\subset \mathbb Z_{\geq 0}.
$$
It is a $\Gamma$-semimodule in the sense that
$$
p\in \Gamma, q\in \Lambda^\mathcal C\Rightarrow p+q\in \Lambda^\mathcal C.
$$
\begin{remark} Note that $\nu_{\mathcal C}(\omega)\geq 1$ for any $\omega\in \Omega^1_{M_0,P_0}$. Anyway, we have the important property that $\Gamma\subset \{0\}\cup \Lambda^{\mathcal C}$. If $\Lambda^{\mathcal C}\cup \{0\}=\Gamma$, we say that $\mathcal C$ is quasi-homogeneous and it is analytically equivalent to the cusp $y^n-x^m=0$.  Otherwise, if $\lambda_1$ is the minimum of  $\Lambda^\mathcal C\setminus \Gamma$, we know that $\lambda_1-n$ is the Zariski invariant, the first nontrivial analytic invariant. This invariant was introduced by Zariski in \cite{Zar1}.
\end{remark}
\begin{remark} Let us note that $\nu_E(\omega)\in \Gamma$, for any $\omega\in \Omega^1_{M_0,P_0}$.
\end{remark}
\subsection{Divisorial order and Differential Values} In view of the definition of the differential values, for any $\omega\in \Omega^1_{M_0,P_0}$ we have that
$
\nu_E(\omega)\leq \nu_{\mathcal C}(\omega)
$.
A useful consequence of this fact is the following one:
\begin{lema}
 \label{lema:resonantbasicforms}
 A basic $1$-form $\omega\in \Omega^1_{M_0,P_0}$ is resonant if and only if $\nu_{\mathcal C}(\omega)>\nu_E(\omega)$.
\end{lema}
\begin{proof} Write $\omega=W+\tilde\omega$, where $W$ is the initial form of $\omega$.
Denote $d=\nu_E(\omega)<nm$. Recall that $\nu_E(\tilde\omega)>\nu_E(W)=d$, we conclude that
$\nu_\mathcal C(\omega)>d$ if and only if $\nu_{\mathcal C}(W)
>d$.
Since $\omega$ is a basic 1-form, we can write
$$
W=x^ay^b\left\{\mu\frac{dx}{x}+\zeta\frac{dy}{y}\right\},\quad an+bm=d.
$$
We have
$$
\phi^*W=(t^n)^a(t^m+t^{m+1}\xi(t))^b\left(n\mu+m\zeta+t\psi(t)\right)\frac{dt}{t}.
$$
The fact that $\nu_\mathcal C(W)>d$ is equivalent to say that $n\mu+m\zeta=0$ and hence it is equivalent to say that $\omega$ is resonant.
\end{proof}

\begin{corolario}
 \label{cor:esonantbasic}
 If $\nu_\mathcal C(\omega)\notin \Gamma$, then $\omega$ is a resonant basic $1$-form.
\end{corolario}
\begin{proof} Since $\nu_{\mathcal C}(\omega)\notin \Gamma$, this differential value is bounded by the conductor $c_\Gamma=(n-1)(m-1)$, hence we have that
$$
\nu_E(\omega)\leq \nu_{\mathcal C}(\omega)<(n-1)(m-1)<nm.
$$
Then $\omega$ is a basic 1-form. Moreover, since $\nu_E(\omega)\in \Gamma$ and $\nu_{\mathcal C}(\omega)\notin\Gamma$, we have that $\nu_E(\omega)<\nu_{\mathcal C}(\omega)$ and we conclude that $\omega$ is a resonant basic 1-form.
\end{proof}

\subsection{Reachability Between Resonant Basic Forms}

Let $\omega,\omega'$ be two $1$-forms $\omega,\omega'\in\Omega^1_{M_0,P_0}$. We say that $\omega'$ is {\em reachable} from $\omega$ if and only if there are nonnegative integer numbers $a,b$  and a constant $\mu\in \mathbb C$ such that
$$
\nu_E(\omega'-\mu x^ay^b\omega)>\nu_E(\omega').
$$
Note that the constant $\mu$ and the pair $(a,b)$ are necessarily unique.

We are interested in the case when $\omega$ and $\omega'$ are basic and resonant. In this situation, the initial parts are respectively given by
$$
W=\mu x^cy^d\left\{m\frac{dx}{dy}-n\frac{dy}{y}\right\},\quad
W'=\mu' x^{c'}y^{d'}\left\{m\frac{dx}{x}-n\frac{dy}{y}\right\}.
$$
Note that $a,a',b,b'\geq 1$ since $\omega,\omega'$ are holomorphic $1$-forms.
We have that $\omega'$ is reachable from $\omega$ if and only if $c'\geq c$ and $d'\geq d$; in this case we have that
$$
\nu_E\left(\omega'-\frac{\mu'}{\mu}x^{c'-c}y^{d'-d}\omega\right)>\nu_E(\omega')=c'n+d'm.
$$
Note also that the minimum divisorial value of a basic and resonant $1$-form is $n+m$ and its initial part is necessarily of the type
$$
\mu xy\left\{m\frac{dx}{x}-n\frac{dy}{y}\right\}=\mu(mydx-nxdy).
$$
If $\omega$ is basic and resonant with $\nu_E(\omega)=n+m$, then any basic and resonant $1$-form is reachable from $\omega$.

\section{Cuspidal Semimodules}
\label{sec:one}
In this section we develop certain features of semimodules over the semigroup $\Gamma$ generated by the Puiseux pair $(n,m)$.
 We consider, unless it is specified, only the singular case $n\geq 2$; in this case the conductor is  $c_\Gamma=(n-1)(m-1)$ and we have the interesting property that any $p\in \Gamma$ with $p<nm$ is written as $p=an+bm$ in a unique way, with $a,b\geq 0$.

We proceed in a self contained way in order to help the reader, several results are true for more general semigroups, but we focus on the cuspidal semigroup $\Gamma$ to shorten the arguments.
\subsection{The Basis of a Semimodule}\label{sec:semi}
A nonempty subset $\Lambda\subset \mathbb Z_{\geq 0}$ is a $\Gamma$-{\em semimodule} if $\Lambda+\Gamma\subset \Lambda$.  We say that $\Lambda$ is {\em normalized} if $0\in \Lambda$, this is equivalent to say that $\Gamma\subset \Lambda$. As for the case of semigroups, the {\em conductor $c_\Lambda$} is defined by
$$
c_\Lambda=\min\{p\in \mathbb Z_{\geq 0};\; \{q\in \mathbb Z; q\geq p \}\subset \Lambda\}.
$$
Note that if $\lambda_{-1}$ is the minimum of $\Lambda$, then we have that $c_\Lambda\leq c_\Gamma+\lambda_{-1}$.

\begin{defn} Let $\Lambda$ be a $\Gamma$-semimodule. A nonempty finite increasing sequence of nonnegative integer numbers
$
\frac{}{}\mathcal B=(\lambda_{-1},\lambda_0,\ldots,\lambda_s)
$
is a {\em basis for} $\Lambda$ if for any $0\leq j\leq s$ we have that
$
\lambda_j\notin \Gamma(\mathcal B_{j-1})
$,
where $\Gamma(\mathcal B_{j-1})=(\lambda_{-1}+\Gamma)\cup(\lambda_0+\Gamma)\cup\cdots\cup (\lambda_{j-1}+\Gamma)$.
\end{defn}
If $\Lambda=\Gamma(\mathcal B)$, we have a chain of semimodules
\begin{equation}
\lambda_{-1}+\Gamma=\Lambda_{-1}\subset \Lambda_{0}\subset\cdots\subset\Lambda_{s}=\Lambda,
\end{equation}
where $\Lambda_j=\Gamma(\mathcal B_j)$. We call {\em decomposition sequence of $\Lambda$} to this chain of semi\-modules. Let us note that
\begin{equation}
\label{eq:unicidadlambdaj}
\lambda_j=\min (\Lambda\setminus \Lambda_{j-1}),\quad 0\leq j\leq s.
\end{equation}
This definitions are justified by next Proposition \ref{prop:basis}
\begin{prop}
\label{prop:basis}
Given a semimodule $\Lambda$, there is a unique basis $\mathcal B$ such that
$
\Lambda=\Gamma(\mathcal B)
$.
\end{prop}
\begin{proof} We start with $\lambda_{-1}=\min \Lambda$. Note that $\Gamma(\lambda_{-1})\subset\Lambda$.  If $\Gamma(\lambda_{-1})=\Lambda$ we stop and we put $s=-1$. If $\Gamma(\lambda_{-1})\ne \Lambda$, we put $\lambda_0=\min(\Lambda\setminus \Gamma(\lambda_{-1}))$.  Note that $\Gamma(\lambda_{-1},\lambda_0)\subset\Lambda$.  We  continue in this way and, since $\lambda_j \not \equiv \lambda_k \mod n$ for $j \neq k$,  after finitely many steps we obtain that
$\Lambda=\Gamma(\lambda_{-1},\lambda_0,\ldots,\lambda_s)$. Let us show the uniqueness of $\mathcal B=(\lambda_{-1},\lambda_{0},\ldots,\lambda_s)$. Assume that $\Lambda=\Gamma(\mathcal B')$, for another $\Gamma$-basis
$
\mathcal B'=(\lambda'_{-1},\lambda'_0,\ldots,\lambda'_{s'})
$.
Note that  $\lambda_{-1}=\min\Lambda=\lambda'_{-1}$.
Assume that $\lambda_j=\lambda'_j$ for any $0\leq j\leq k-1$. In view of Equation \eqref{eq:unicidadlambdaj} we have that
$
\lambda_k=\lambda'_k=\min\left( \Lambda\setminus \Gamma(\mathcal B_{k-1})\right)=
\min\left( \Lambda\setminus \Gamma(\mathcal B'_{k-1})\right)
$.
This ends the proof.
\end{proof}

We say that $\mathcal B=(\lambda_{-1},\lambda_0,\ldots,\lambda_s)$ is {\em the basis of $\Lambda=\Gamma(\mathcal B)$} and that $s$ is the {\em length} of $\Lambda$.

Consider a semimodule $\Lambda=\Gamma(\mathcal B)$, an element $\lambda\in \mathbb Z_{\geq 0}$ is said to be {\em $\Lambda$-independent} if and only if $\lambda\notin \Lambda$ and $\lambda>\lambda_s$, where $\lambda_s$ is the last element in the basis $\mathcal B$. In this case we obtain a basis $\mathcal B(\lambda)$, just by adding $\lambda$ to $\mathcal B$ as being the last element. The new semimodule is denoted $\Lambda(\lambda)$, thus we have
$
\Lambda(\lambda)=\Lambda\cup (\lambda+\Gamma)=\Gamma(\mathcal B(\lambda))
$.

Given a  semimodule $\Lambda=\Gamma(\lambda_{-1},\lambda_0,\ldots,\lambda_s)$, we define the {\em axes} $u_i=u_i(\Lambda)$ by
\begin{equation}
u_0=\lambda_{-1};\quad
u_i=\min \left(\Lambda_{i-2}\cap (\lambda_{i-1}+\Gamma)\right), \; 1\leq i\leq s+1.
\end{equation}
Note that
$
u_i(\Lambda_j)=u_i(\Lambda)$, for $0\leq i\leq j+1\leq s+1
$.

\begin{defn}\label{def:insemi:increasing}
A semimodule $\Lambda=\Gamma(\lambda_{-1},\lambda_0,\ldots,\lambda_s)$ is {\em increasing} if and only if
$
\lambda_i>u_i
$ for any $i=0,1,\ldots,s$.
\end{defn}

\begin{remark}\label{re:insemi:decomposition}
If $\Lambda$ is an increasing semimodule, each element $\Lambda_i$ of the decomposition sequence is also  an  increasing semimodule. Moreover, if $\lambda'$ is a $\Lambda$-independent value with $\lambda'>u_{s+1}$, then $\Lambda(\lambda')$ is also an increasing $\Gamma$-semimodule.
\end{remark}

Given a semimodule $\Lambda=\Gamma(\lambda_{-1},\lambda_0,\ldots,\lambda_s)$, the semimodule $\widetilde{\Lambda}=\Lambda-\lambda_{-1}$ is called the {\em normalization  of $\Lambda$}. Next features allow to deduce properties of $\Lambda$ from properties of its normalization:
\begin{enumerate}
\item  The basis of $\widetilde{\Lambda}$ is $(0,\lambda_0-\lambda_{-1},\ldots,\lambda_s-\lambda_{-1})$.
\item $\widetilde{\Lambda}_i=\Lambda_i-\lambda_{-1}$, for $i=-1,0,\ldots,s$.
\item $
u_i(\widetilde{\Lambda})=u_i(\Lambda)-\lambda_{-1},\quad i=0,1,\ldots,s+1
$.
\item $ c_{\widetilde{\Lambda}}=c_\Lambda-\lambda_{-1}$.
\item $\Lambda$ is increasing if and only if $\widetilde{\Lambda}$ is increasing.
\end{enumerate}

\subsection{Axes and conductor} We precise the expressions of the axes and we bound them by the conductors.
\begin{lema}
 \label{lema:cotaconductor}
 Consider a semimodule $\Lambda$ of length $s$ and two indices $0\leq k<i\leq s$. Then we have that $u_{i+1}<c_{\Lambda_k}+n$.
\end{lema}
\begin{proof}
   We can assume that $i=s$, $k=s-1$. Note that $\lambda_{s}<c_{\Lambda_{s-1}}$, since $\lambda_s\notin \Lambda_{s-1}$. Then, there is a unique
 $\alpha\in \mathbb Z_{> 0}$ such that
$
0\leq \lambda_{s}-c_{\Lambda_{s-1}}+\alpha n<n
$.
We have that $\lambda_{s}+\alpha n\in \lambda_{s}+\Gamma$ and $\lambda_{s}+\alpha n\geq c_{\Lambda_{s-1}}$. We obtain
$$
\lambda_{s}+\alpha n\in \Lambda_{s-1}\cap (\lambda_{s}+\Gamma).
$$
We deduce that $u_{s+1}\leq \lambda_{s}+\alpha n<c_{\Lambda_{s-1}}+n$.
\end{proof}

\begin{corolario}
 \label{cor:cotaui2}
 Consider a semimodule $\Lambda$ of the form $\Lambda=\Gamma(n,m,\lambda_1,\ldots,\lambda_s)$. Then $u_{i+1}<nm$, for any $0\leq i\leq s$.
\end{corolario}
\begin{proof} By Lemma~\ref{lema:cotaconductor}, we have that $u_{i+1}\leq c_{\Lambda_0}+n$, but in this situation, we have that $\Lambda_0\cup\{0\}=\Gamma$ and thus
$$
c_{\Lambda_0}=c_\Gamma=(n-1)(m-1).
$$
Hence $u_{i+1}\leq c_{\Lambda_0}+n=(n-1)(m-1)+n<nm$.
\end{proof}
\begin{lema}\label{re:insemi:ui} Consider $\Lambda=\Gamma(\lambda_{-1},\lambda_0,\ldots,\lambda_s)$. There is a unique index $k$ with $-1\leq k\leq s-1$ such that $u_{s+1}\in \lambda_k+\Gamma$ and there are unique expressions \begin{eqnarray}
\label{eq:1}
u_{s+1}&=&\lambda_{s}+na+mb,\quad  a,b\in \mathbb Z_{\geq 0}
\\
\label{eq:2}
u_{s+1}&=&\lambda_k+n c+m d,\quad
 c,d\in \mathbb Z_{\geq 0}.
\end{eqnarray}
In these expressions, we have $ac=bd=ab=cd=0$ and $(a,b)\ne(0,0)\ne(c,d)$.
\end{lema}

\begin{proof}
The existence of the expressions \eqref{eq:1} and \eqref{eq:2} is given by the definition of $u_{s+1}$ as $u_{s+1}=\min (\Lambda_{s-1}\cap (\lambda_s+\Gamma))$. By the minimality of $u_{s+1}$ and the fact that $u_{s+1}\ne\lambda_{s}$ and $u_{s+1}\ne\lambda_{k}$, we deduce the properties $ac=bd=0$ and $(a,b)\ne(0,0)\ne(c,d)$. Moreover, if $ab\ne 0$ we should have that $c=d=0$ which is not possible; in the same way we see that $cd=0$.

Let us show the uniqueness of the index $k$. Assume that there are two indices
 $-1\leq k< k'\leq s-1$ with
 $
 u_{s+1}=\lambda_k+c n+d m= \lambda_{k'}+c'n+d'm
 $.
 Take the case when $a\ne0$, then we have that $c=c'=0$ and we can write
 $$
 \lambda_{k'}=\lambda_k+m(d-d')\in \lambda_k+\Gamma,
 $$
 this is a contradiction. Same argument if $b\ne 0$.

If we normalize $\Lambda$, we have that $u_{s+1}-\lambda_{s}=\tilde u_{s+1}-\tilde \lambda_{s}$. By Lemma \ref{lema:cotaconductor} we have
 $$\tilde u_{s+1}-\tilde\lambda_{s}\leq
 \tilde u_{s+1}< c_{\widetilde{\Lambda}_{s-1}}+n\leq c_\Gamma+n=(n-1)(m-1)+n<nm.
 $$
Hence, we have that $u_{s+1}-\lambda_{s}$ (and with the same argument $u_{s+1}-\lambda_k$) are strictly smaller than $nm$. Thus, the expression of these elements of $\Gamma$ as a linear combination of $n,m$ with non-negative coefficients is unique.
\end{proof}

\subsection{The limits} Consider a semimodule $\Lambda=\Gamma(\lambda_{-1},\lambda_0,\ldots,\lambda_s)$ with $s\geq 0$.
The {\em first and second limits $\ell_1$ and $\ell_2$ of $\Lambda$} are defined by
\begin{eqnarray}
\ell_1&=&\min\{p\,;\;np+\lambda_s\in \Lambda_{s-1}\}.\\
\ell_2&=&\min\{q\,;\; mq+\lambda_s\in \Lambda_{s-1}\}.
\end{eqnarray}
\begin{remark}
\label{rk:expresionesui}
We have that $\ell_1\ell_2\ne 0$ and
\begin{equation}
\label{eq:uesemasuno}
u_{s+1}=\min\{\ell_1n+\lambda_s,\; \ell_2m+\lambda_s\}.
\end{equation}
Indeed, by Lemma \ref{re:insemi:ui}, we have either $u_{s+1}=an+\lambda_s$ or $u_{s+1}=bm+\lambda_{s}$; if $u_{s+1}=an+\lambda_s$, by minimality we have that $a=\ell_1$, in the same way, if $u_{s+1}=bm+\lambda_s$ we have that $b=\ell_2$. Moreover, there is a unique index $k$ with $-1\leq k<s$ such that
\begin{enumerate}
\item  If $u_{s+1}=\ell_1n+\lambda_s$, then $u_{s+1}=\lambda_k+bm$.
\item  If $u_{s+1}=\ell_2 m+\lambda_s$, then $u_{s+1}=\lambda_k+an$.
\end{enumerate}
\end{remark}
\begin{lema}
 \label{lema:limits}
 If $an+bm+\lambda_s\in \Lambda_{s-1}$, then either $a\geq \ell_1$ or $b\geq \ell_2$.
\end{lema}
\begin{proof} Let us write $an+bm+\lambda_s=cn+dm+\lambda_{j}$ for a certain $j\leq s-1$. If $ac\ne 0$, we find
$$(a-1)n+bm+\lambda_s=(c-1)n+dm+\lambda_j\in \Lambda_{s-1}.$$
 Repeating the argument and working in a similar way with the coefficients $b,d$, we find an element
$$
\tilde a n+\tilde b m+\lambda_{s}= \tilde c n+\tilde d m+\lambda_j
$$
such that $a\geq \tilde a$ and $b\geq \tilde b$, with the property that $\tilde a \tilde c=0$ and $\tilde b\tilde d=0$. Moreover, we have that $(\tilde c,\tilde d)\ne (0,0)$, since otherwise $\lambda_s\leq \lambda_{s-1}$. Suppose that $\tilde c\ne 0$, then $\tilde a=0$ and
$$
\tilde b m+\lambda_{s}= \tilde c n+\tilde d m+\lambda_j\in \Lambda_{s-1}.
$$
By the minimality property of $\ell_2$, we have that $\tilde b\geq\ell_2$ and then $b\geq \ell_2$. In a similar way, we show that if $\tilde d\ne 0$ we have that $a\geq\ell_1$.
\end{proof}
Let us note that the limits of the normalization $\tilde\Lambda$ are the same ones as for $\Lambda$.
\begin{example}\label{ex:5-11} Consider the semigroup $\Gamma=\langle 5,11\rangle$ and the $\Gamma$-semimodule \linebreak $\Lambda=\Gamma(5,11,17,23,29)$. Let us compute the axes and the limits for this semimodule. Note that $s=3$. We have $u_0=\lambda_{-1}=n=5$, and $\lambda_0=m=11$. In order to compute the limits of $\Lambda_{0}$ we have to find the minimal non-negative integers
\begin{itemize}
  \item $\ell_1^1$   such that $11+5 \ell_1^1=5 +11b$, hence $\ell_1^1=b=1$;
  \item $\ell_2^1$   such that $11+11 \ell_2^1=5+5a$ and we obtain $\ell_2^1=4$ and $a=10$.
\end{itemize}
Hence, $\lambda_0+n \ell_1^1=16$ and $\lambda_0+m \ell_2^1=55$ and $u_1=\min\{16,55\}=16$. Now, let us compute the limits of $\Lambda_{1}$ where $\lambda_1=17$. We search $\ell_1^2$ and $\ell_2^2$ minimal such that
\begin{itemize}
  \item $17+5 \ell_1^2=11+11b$, hence $\ell_1^2=b=1$ and $\lambda_1+n\ell_1^2=22$;
  \item $17+11\ell_2^2=5 + 5a$ and we have that $\ell_2^2=3$ and $a=9$, then $\lambda_1+m\ell_2^2=50$.
\end{itemize}
We get that $u_2=22$. In a similar way, taking into account that $\lambda_2=23$ and $\lambda_3=29$, we get that $u_3=28$ and $u_4=34$. Since $u_i < \lambda_i$, for $i=0,1,2,3$,
 we obtain that the semimodule $\Lambda$ is increasing.
\end{example}

\section{Standard Bases}\label{sec:differe}
From now on, we fix a cusp $\mathcal C$ in $\operatorname{Cusps}(E)$ and we consider the semimodule $\Lambda$ of differential values of $\mathcal C$:
\begin{equation}
\Lambda=\Lambda^\mathcal C=\{\nu_\mathcal C(\omega);\; \omega \in \Omega^1_{M_0,P_0}\}\setminus\{\infty\}.
\end{equation}
We recall that $\Gamma\setminus \{0\}\subset \Lambda$.
\begin{lema} If $(\lambda_{-1},\lambda_0,\lambda_1,\ldots,\lambda_s)$ is the basis of $\Lambda^\mathcal C$, then $\lambda_{-1}=n$ and $\lambda_0=m$.
\end{lema}
\begin{proof} Let $(x,y)=(t^n, t^m\xi(t))$ be a Puiseux parametrization of $\mathcal C$, where $\xi(0)\ne 0$. Recall that $\nu_\mathcal C(adx+bdy)$ is the order in $t$ of the expression
\begin{equation}
\label{eq:valordiferencial}
nt^{n}a(t^n, t^m\xi(t))+
t^{m}\xi(t)b(t^n, t^m\xi(t))\{m+t\xi'(t)/\xi(t)\}.
\end{equation}
We see that this order is $\geq n$ and that $\nu_{\mathcal C}(dx)=n$. Hence $n=\lambda_{-1}$. Moreover, the terms in Equation \eqref{eq:valordiferencial} of degree $<m$ come only from the first part  $nt^{n}a(t^n, t^m\xi(t))$ of the sum, so, they are values in $\Gamma$. Since $m=\nu_\mathcal C(dy)$, we conclude that $\lambda_0=m$.
\end{proof}
\begin{defn} Write $\Lambda^\mathcal C=\Gamma(n,m,\lambda_1,\ldots,\lambda_s)$. A {\em standard basis for $\mathcal C$} is a list of $1$-forms
$
\mathcal G=(\omega_{-1},\omega_0,\omega_1,\ldots,\omega_s)
$
such that $\nu_{\mathcal C}(\omega_i)=\lambda_i$, for $i=-1,0,1,\ldots,s$.
\end{defn}
\begin{remark} There is at least one standard basis, by definition of the semimodule of differential values. The standard bases are not in general unique. For instance, we have that
$$
\nu_{\mathcal C}(h dx)=n,\quad \nu_{\mathcal C}(h dy)=m,\quad h(0)\ne 0.
$$
On the other hand, for $i=1,2,\ldots,s$, we have that $\nu_\mathcal C(\omega_i)=\lambda_i\notin \Gamma$, then, in view of Corollary  \ref{cor:esonantbasic}, the 1-form $\omega_i$ is basic resonant.
\end{remark}
\begin{remark} Let $\mathcal G=(\omega_{-1},\omega_0,\omega_1,\ldots,\omega_s)$ be a standard basis for $\mathcal C$. Then $\omega_{-1},\omega_0$ have the form
$$
\omega_{-1}=h dx+gdy, \;h(0)\ne 0;\quad
\omega_0= f dx+\psi dy,\;  \psi(0)\ne 0,\; \nu_\mathcal C(fdx)>m.
$$
Thus, we can write any differential $1$-form $\omega$ in a unique way as
$\omega=a\omega_{-1}+b\omega_{0}$. Anyway, we are mainly interested in the study of the $1$-forms $\omega_i$, for $1\leq i\leq s$. From the standard basis $\mathcal G$ we can obtain a new one ``adapted to the coordinates'' given by
$$
\mathcal G=(dx,dy,\omega_1,\ldots,\omega_s).
$$
Just for simplifying the presentation of the computations, we will consider only this kind of standard bases.
\end{remark}

\subsection{The Zariski Invariant}

In this subsection, we deal with properties of divisorial orders  and  differential values around the element $\lambda_1$, where
$$\Lambda^\mathcal C=\Gamma(n,m,\lambda_1,\ldots,\lambda_s).
$$ This is the first step for a general result.
Anyway, let us recall that $\lambda_1-n$ is the classical Zariski invariant. Let us cite the work of O. Gómez \cite{Oziel} that essentially contains several of the results in this Section.

\begin{prop}
\label{prop:ordmonzariski} We have the following properties:
\begin{enumerate}
\item If $s=0$, then $\infty= \sup\{\nu_{\mathcal C}(\omega);\; \omega\in \Omega^1_{M_0,P_0},\;\nu_E(\omega)=n+m\}$.
\item If $s\geq 1$, then $\lambda_1= \sup\{\nu_{\mathcal C}(\omega);\; \omega\in \Omega^1_{M_0,P_0},\;\nu_E(\omega)=n+m\}$.
\end{enumerate}
\end{prop}
\begin{proof} Assume that $s=0$ and hence $\Lambda^\mathcal C=\Gamma\setminus\{0\}$. Let us consider the $1$-form $\eta=mydx-nxdy$. We have that $\nu_{\mathcal C}(\eta)>n+m=\nu_E(\eta)$. Moreover, since $s=0$ we have that $\nu_{\mathcal C}(\eta)\in \Gamma$; then there is a monomial function $f$ such that $$\nu_E(df)=\nu_\mathcal C(df)= \nu_{\mathcal C}(\eta)>n+m.
$$
 In particular, there is a constant $\mu\ne 0$ such that
$
\nu_\mathcal C(\eta-\mu df)>\nu_{\mathcal C}(\eta)
$. Write $\eta^1=\eta-\mu df$; we have that $\nu_E(\eta^1)=\nu_E(\eta)=n+m$ and $\nu_\mathcal C(\eta^1)>\nu_{\mathcal C}(\eta)$. We repeat the argument with $\eta^1$ and in this way we obtain $1$-forms $\eta^k$ with $\nu_E(\eta^k)=n+m$ and $\nu_{\mathcal C}(\eta^k)\geq n+m+1+k$. This proves the first statement.

Assume now that $s\geq 1$. Let us first show that
 $$
 \lambda_1\leq \sup\{\nu_{\mathcal C}(\omega);\; \omega\in \Omega^1_{M_0,P_0},\;\nu_E(\omega)=n+m\}.
 $$

If $\nu_{\mathcal C}(\eta)\notin \Gamma$, we have that
$\nu_{\mathcal C}(\eta)\geq \lambda_1$ since $\lambda_1$ is the minimum of the differential values not in $\Gamma$, then we are done.
 Assume that $\nu_{\mathcal C}(\eta)\in \Gamma$ and hence $$\nu_{\mathcal C}(\eta)=an+bm>n+m.$$ Taking the function $f=x^ay^b$, up to multiply $df$ by a constant $c_1$ we obtain that
$$
\nu_{\mathcal C}(\eta_1)>\nu_{\mathcal C}(\eta)=an+bm, \quad \eta_1=\eta-c_1df.
$$
Note that $\nu_E(\eta_1)=n+m$, since
$
\nu_E(df)=an+bm>n+m.
$
 We restart with $\eta_1$ instead of $\eta$, noting that $\nu_{\mathcal C}(\eta)<\nu_{\mathcal C}(\eta_1)$. Repeating finitely many times this procedure,  we obtain a new $1$-form $\tilde\eta=\eta-d\tilde f$ such that $\nu_E(\tilde \eta)=n+m$ and either $\nu_{\mathcal C}(\tilde\eta)\geq c_\Gamma=(n-1)(m-1)$ or $\nu_{\mathcal C}(\tilde\eta)\notin \Gamma$, in both cases we have that $\nu_{\mathcal C}(\tilde\eta)\geq\lambda_1$ and we are done.

 It remains to show that $
 \lambda_1\geq \sup\{\nu_{\mathcal C}(\omega);\; \omega\in \Omega^1_{M_0,P_0},\;\nu_E(\omega)=n+m\}
 $.
  Let us consider $\omega_1$ such that $\nu_{\mathcal C}(\omega_1)=\lambda_1$ and let us show that it is not possible to have  $\tilde\omega$ such that $\nu_E(\tilde\omega)=n+m$ and $\nu_{\mathcal C}(\tilde \omega)>\nu_{\mathcal C}(\omega_1)$. In this situation, both $\omega_1$ and $\tilde\omega$ are basic resonant. We know that $\omega_1$ is reachable from $\tilde\omega$ and thus there is a constant $\mu$ and $a,b\geq 0$ such that
 $$
 \nu_E(\omega_1^1)>\nu_E(\omega_1), \quad
 \omega^1_1= \omega_1-\mu x^ay^b\tilde\omega.
 $$
 We have that $\nu_{\mathcal C}(\omega_1^1)=\nu_{\mathcal C}(\omega_1)=\lambda_1$. We restart with the pair $\omega_1^1,\tilde\omega$; in this way, we obtain an infinite sequence of $1$-forms
 $
 \omega_1,\omega_1^1,\omega_1^2,\ldots
 $
 with strictly increasing divisorial orders.
Up to a finite number of steps, we find an index $k$ such that $\nu_E(\omega_1^k)>\lambda_1=\nu_{\mathcal C}(\omega_1^k)$.
 This contradicts with the fact $\nu_{\mathcal C}(\omega_1^k)\geq\nu_E(\omega_1^k)$.
\end{proof}

\begin{corolario} Any $1$-form $w\in \Omega^1_{M_0,P_0}$ such that $\nu_E(\omega)=n+m$ and $\nu_{\mathcal C}(\omega)\notin \Gamma$ satisfies that $\nu_{\mathcal C}(\omega)=\lambda_1$.
\end{corolario}

\begin{proof} In view of the previous result, we have that $\nu_{\mathcal C}(\omega)\leq \lambda_1$. Since $\nu_{\mathcal C}(\omega)\notin\Gamma$, we also have that
$\nu_{\mathcal C}(\omega)\geq \lambda_1$.
\end{proof}

\begin{corolario}
 \label{cor:ordlambdauno}Any $1$-form $w\in \Omega^1_{M_0,P_0}$ such that $\nu_\mathcal C(\omega)=\lambda_1$ satisfies that $\nu_E(\omega)=n+m$.
\end{corolario}
\begin{proof} Take $\omega_1$ such that $\nu_\mathcal C(\omega_1)=\lambda_1$ and $\nu_E(\omega_1)=n+m$. Assume that $$\nu_E(\omega)>n+m$$ in order to obtain a contradiction.
 Since $\lambda_1\notin\Gamma$, both $\omega$ and $\omega_1$ are basic resonant and $\omega$ is reachable from $\omega_1$. Then there is a function $f$ with $\nu_{\mathcal C}(f)>0$ such that
$$
\nu_E(\omega-f\omega_1)>\nu_E(\omega).
$$
Put $\omega^1=\omega-f\omega_1$, since $\nu_{\mathcal C}(f\omega_1)>\lambda_1$, we have that $\nu_{\mathcal C}(\omega^1)=\lambda_1$. We restart with the pair $\omega^1,\omega$. After finitely many repetitions we find $\omega^k$ with $\nu_\mathcal C(\omega^k)=\lambda_1$ and
 $\nu_E(\omega^k)>\lambda_1$, contradiction.
\end{proof}

The following two lemmas are necessary steps in order to prove an inductive version of Proposition \ref{prop:ordmonzariski} valid for all indices $i=1,2,\ldots,s$:

\begin{lema}\label{le:bs:w14} Assume that $s\geq 1$ and take $\omega_1$ such that $\nu_{\mathcal C}(\omega_1)=\lambda_1$. Consider an integer number
$
k=na+mb+\lambda_1\in \lambda_1+\Gamma
$.
 The following statements are equivalent:
 \begin{enumerate}
 \item $k\notin \Gamma$.
 \item  $\nu_{E}(\omega)\leq\nu_{E}(x^ay^b\omega_1)$, for
any $\omega\in \Omega^1_{M_0,P_0}$  such that $\nu_\mathcal C(\omega)=k$.
 \end{enumerate}
\end{lema}

\begin{proof}
Note that $k=\nu_\mathcal C(x^ay^b\omega_1)>(a+1)n+(b+1)m=\nu_E(x^ay^b\omega_1)$.

 Assume that $k\in \Gamma$, then $k=na'+mb'>\nu_{E}(x^ay^b\omega_1)$. Taking $\omega=d(x^{a'}y^{b'})$, we have $\nu_{\mathcal C}(\omega)=\nu_E(\omega)=k>\nu_E(x^ay^b\omega_1)$.

Now assume that $k\notin\Gamma$. Let us reason by contradiction assuming that there is $\omega$ with $\nu_{\mathcal C}(\omega)=k$ with $\nu_{E}(\omega)> \nu_{E}(x^ay^b\omega_1)$.
We have that $\omega$ is basic resonant, since $\nu_\mathcal C(\omega)\notin \Gamma$. Then $\omega$ is reachable from $\omega_1$. Then there is $a',b'\geq 0$ and a constant $\mu$ such that $\nu_E(x^{a'}y^{b'}\omega_1)=\nu_E(\omega)$ and
$$
\nu_E(\omega-cx^{a'}y^{b'}\omega_1)>\nu_E(\omega)>\nu_E(x^{a}y^{b}\omega_1).
$$
Since $na'+mb'>na+mb$,  we have that  $\nu_{\mathcal C}(x^{a'}y^{b'}\omega_1)>k$ and hence $\nu_{\mathcal C}(\omega^1)=k$, where $\omega^1=\omega-cx^{a'}y^{b'}\omega_1$. Repeating the procedure with the pair $\omega^1,\omega_1$, we obtain a sequence
$$
\omega,\omega^1,\omega^2,\ldots
$$
with strictly increasing divisorial order and such that $\nu_\mathcal C(\omega^j)=k$ for any $j$. This is a contradiction.
\end{proof}

\begin{lema}\label{le:bs:w15} Take $\omega_1$ with $\nu_{\mathcal C}(\omega_1)=\lambda_1$.
Let $\omega\in \Omega^1_{M_0,P_0}$ be a 1-form such that $\nu_\mathcal C(\omega)=\lambda\notin \Gamma$. There are unique $a,b\geq 0$ such that $\nu_{E}(\omega)=\nu_{E}(x^ay^b\omega_1)$. Moreover, we have that $\lambda\geq na+mb+\lambda_1$.
\end{lema}

\begin{proof} Note that $\omega$ is basic resonant and thus the existence and uniqueness of $a,b$ is assured. Moreover, if $\lambda<na+mb+\lambda_1$, we can find a constant $\mu$ such that
$$\nu_{E}(\omega-\mu x^ay^b\omega_1)> \nu_{E}(x^ay^b\omega_1)
$$
and
$\nu_\mathcal C(\omega-\mu x^ay^b\omega_1)=\lambda$. Put $\omega^1=\omega-\mu x^ay^b\omega_1$, we have that $\nu_{\mathcal C}(\omega^1)=\lambda\notin \Gamma$. As before, we have that
$$
\nu_E(\omega^1)=\nu_E(x^{a_1}y^{b_1}\omega_1),\quad a_1n+b_1m>an+bm
$$
and thus $\lambda<a_1n+b_1m+\lambda_1$.
We repeat the process with the pair $\omega^1,\omega_1$, where $\omega^1=\omega-\mu x^ay^b\omega_1$  in order to have a sequence $\omega,\omega^1,\omega^2,\ldots $ with strictly increasing divisorial orders and such that $\nu_{\mathcal C}(\omega^j)=\lambda$ for any $j$. This is a contradiction.
\end{proof}

\subsection{Critical Divisorial Orders} Recall that we are considering a cusp $\mathcal C$ in $\operatorname{Cusps}(E)$, whose semimodule of differential values is
$$
\Lambda^{\mathcal{C}}=\Gamma(n,m,\lambda_1,\ldots,\lambda_s).
$$
The {\em critical divisorial orders  $t_i$}, for $i=-1,0,\ldots,s+1$ are defined as follows:
\begin{itemize}
\item We put $t_{-1}=n$ and  $t_0=m$.
\item For $1\leq i\leq s+1$, we put $t_i=t_{i-1}+u_i-\lambda_{i-1}$.
\end{itemize}
Let us note that $t_1=m+(n+m)-m=n+m$.

\begin{lema}
\label{lema:lambdamenoslambda} Consider the semimodule $\Lambda=\Gamma(n,m,\lambda_1,\ldots,\lambda_s)$ and take an index $1\leq i\leq s$.
If $\lambda_\ell>u_\ell$, for any $0\leq \ell\leq i$, we have that
\begin{equation}
\lambda_j-\lambda_k> t_j-t_k,\quad -1\leq k< j\leq i.
\end{equation}
\end{lema}
\begin{proof} We have that $\lambda_{j}-\lambda_{j-1}>t_{j}-t_{j-1}$ if and only if
$$
t_j=t_{j-1}+u_j-\lambda_{j-1}>t_{j}+u_j-\lambda_j,
$$ which is true, since $u_j-\lambda_j<0$. Noting that
$$
\lambda_j-\lambda_k=\sum_{\ell=k}^{j-1}(\lambda_{\ell+1}-\lambda_\ell)>
\sum_{\ell=k}^{j-1}(t_{\ell+1}-t_\ell)=t_j-t_k,
$$
The proof is ended. \end{proof}

\begin{lema}
\label{lema:lambdamenoslambda2} If the semimodule $\Lambda=\Gamma(n,m,\lambda_1,\ldots,\lambda_s)$ is increasing, we have that $t_i<nm$, for any $i=-1,0,1,\ldots,s+1$.
\end{lema}
\begin{proof} If $i\in\{-1,0\}$ we have that $t_{-1}=n$, $t_0=m$ and we are done. Assume that $1\leq i\leq s$, we have that
$$t_i-n=
t_i-t_{-1}=\sum_{\ell=0}^i(t_\ell-t_{\ell-1})\leq \sum_{\ell=0}^i(\lambda_\ell-\lambda_{\ell-1})=\lambda_{i}-\lambda_{-1}=\lambda_i-n.
$$
Then $t_i\leq \lambda_i<c_\Gamma<nm$. Consider the case $i=s+1$. We have that $$t_{s+1}=t_s+u_{s+1}-\lambda_s=u_{s+1}+(t_s-\lambda_s)\leq u_{s+1}<c_\Gamma+n<nm.
$$
See Lemma \ref{lema:cotaconductor}.
\end{proof}

\begin{remark}
As a consequence of Lemma \ref{lema:lambdamenoslambda2} we have that any $1$-form $\omega$ such that $\nu_E(\omega)=t_i$ is a basic 1-form; moreover, if $t_i=\nu_E(\omega)<\nu_\mathcal C(\omega)$, then it is basic and resonant.
\end{remark}

The critical divisorial orders are the divisorial orders of the elements of a standard basis, in view of the following

\begin{theorem}\label{teo:bs:delorme}
For each $1\leq i\leq s$ we have the following statements
\begin{enumerate}
\item $\lambda_i=\sup\{\nu_\mathcal C(\omega):\nu_{E}(\omega)=t_i\}$.
\item If $\nu_\mathcal C(\omega)=\lambda_i$, then $\nu_{E}(\omega)=t_i$.
\item For each $1$-form $\omega$ with $\nu_\mathcal C(\omega)\notin \Lambda_{i-1}$, there is a unique pair $a,b\geq 0$ such that $\nu_{E}(\omega)=\nu_{E}(x^ay^b\omega_i)$. Moreover, we have that $\nu_\mathcal C(\omega)\geq \lambda_i+na+mb$.
\item We have that $\lambda_i>u_i$.
\item Let $k=\lambda_i+na+mb$, then $k\notin \Lambda_{i-1}$ if and only if for all $\omega$ such that $\nu_\mathcal C(\omega)=k$ we have that $\nu_{E}(\omega)\leq \nu_{E}(x^ay^b\omega_i)$.
\end{enumerate}
In particular, the semimodules $\Lambda_i$ are increasing, for $i=1,2,\ldots,s$.
\end{theorem}
A proof of this Theorem \ref{teo:bs:delorme} is given in Appendix \ref{Apendice-Semimodulo}.
\begin{remark} Note that if
$\mathcal B=(\omega_{-1}=dx,\omega_0=dy,\omega_1,\omega_2,\ldots,\omega_s)$
is a standard basis, Theorem \ref{teo:bs:delorme} says that $\nu_E(\omega_i)=t_i$ for any
$i=-1,0,1,\ldots,s$ and that $\omega_{j+1}$ is reachable from $\omega_j$, for any $1\leq j\leq s-1$. That is, the initial parts of the $1$-forms $\omega_i$ are given by
$$
W_i=\mu_ix^{a_i}y^{b_i}\left\{m\frac{dx}{x}-n\frac{dy}{y}\right\},\quad \mu_i\ne 0,
$$
where $(a_1,b_1)=(1,1)$ and
$$
1=a_1\leq a_2\leq\cdots \leq a_s,\quad 1=b_1\leq b_2\leq\cdots \leq b_s.
$$
Moreover, we have that
$
na_i+mb_i=t_i$, for $ i=1,2,\ldots,s
$.
\end{remark}
\begin{remark} Note that $\mathcal B=(\omega_{-1}=dx,\omega_0=dy,\omega_1,\omega_2,\ldots,\omega_s)$
is a standard basis if and only if $\nu_E(\omega_i)=t_i$ and $\nu_\mathcal C(\omega_i)\notin \Lambda_{i-1}$, for any $i=1,2,\ldots,s$. Moreover, Theorem \ref{teo:bs:delorme} justifies an algorithm of construction of a standard basis as follows:
\begin{quote}\em  Assume we have obtained $\omega_j$, for $j=-1,0,1,\ldots,s'$. We can produce the axis $u_{s'+1}$ and the critical divisorial order $t_{s'+1}=t_{s'}+u_{s'+1}-\lambda_{s'}$. There is an expression
$t_{s'+1}=an+bm$. We consider the $1$-form
$$
\omega_{s'+1}^0=x^ay^b\left\{m\frac{dx}{x}-n\frac{dy}{y}\right\}.
$$
If $\nu_\mathcal C(\omega_{s'+1}^0)\notin \Lambda_{s'}$, we know that $\lambda_{s'+1}=\nu_\mathcal C(\omega_{s'+1}^0)$ and $s\geq s'+1$. If $\nu_\mathcal C(\omega_{s'+1}^0)\in \Lambda_{s'}$, there is $j\leq s'$ and $c,d\geq 0$ such that
$$
\nu_\mathcal C(x^cy^d\omega_j)= \nu_\mathcal C(\omega_{s'+1}^0).
$$
We take a constant $\mu$ such that $\nu_\mathcal C (\omega_{s'+1}^0-\mu x^cy^d\omega_j)>\nu_\mathcal C(\omega_{s'+1}^0)$. Put $\omega_{s'+1}^1= \omega_{s'+1}^0-\mu x^cy^d\omega_j$. We have that $\nu_E(\omega_{s'+1}^1)=t_{s'+1}$. We repeat the procedure with $\omega_{s'+1}^1$. After finitely many steps we get that either $\nu_\mathcal C(\omega_{s'+1}^k)\notin \Lambda_{s'}$ or $\nu_\mathcal C(\omega_{s'+1}^k)\geq c_\Gamma$. In the first case, we put $\lambda_{s'+1}=\nu_\mathcal C(\omega_{s'+1}^k)$,  in the second case we know that $s=s'$.
\end{quote}
\end{remark}

\begin{example}\label{ex:5-11-v2}
Consider the semigroup $\Gamma=\langle 5,11\rangle$ and the $\Gamma$-semimodule \linebreak $\Lambda=\Gamma(5,11,17,23,29)$ as in Example~\ref{ex:5-11}. The computation of the critical divisorial orders $t_i$, $i=-1,0,1,2,3$, gives
$$t_{-1}=5, \quad t_0=11, \quad t_1=16, \quad t_2=21,  \quad t_3=26, \quad t_4=31.$$
Since the semimodule $\Lambda$ is increasing, by \cite{patricio} there exists a curve whose semimodule is $\Lambda$. In particular,   the curve $\mathcal{C}$  given by the Puiseux parametrization \linebreak $\phi(t)=(t^5,t^{11}+t^{12}+t^{13})$  has semigroup  $\Gamma$ and semimodule of differential values $\Lambda$. An  extended standard basis is given by $\{\omega_{-1},\omega_0,\omega_1,\omega_2,\omega_3,\omega_4\}$ where
$$\omega_{-1}=dx, \quad \omega_0=dy, \quad \omega_1=5x dy-11ydx, \quad \omega_2=11 x \omega_1 -5y dy, \quad \omega_3=x \omega_2+y\omega_1$$
and $\omega_4=x \omega_3 - 33 y \omega_2-1199x^6 dx -2035x^4 y dx -407 x^4 \omega_1 -1595 x^3 \omega_2 + \cdots$. The reader can check  that $\nu_E(\omega_i)=t_i$ as stated in Theorem~\ref{teo:bs:delorme}.
\end{example}

\section{Extended Standard Basis and Analytic Semiroots}
As in previous sections, we consider a cusp $\mathcal C$ in $\operatorname{Cusps}(E)$, whose semimodule of differential values is
$$
\Lambda=\Gamma(n,m,\lambda_1,\ldots,\lambda_s).
$$
Let us recall that the axis $u_{s+1}=\min (\Lambda_{s-1}\cap (\lambda_s+\Gamma))$ is well defined and we have also a well defined critical divisorial order
$$
t_{s+1}=t_s+u_{s+1}-\lambda_s.
$$
 Let us also remark that $t_\ell<nm$, for $0\leq \ell\leq s+1$, in view of Lemma \ref{lema:lambdamenoslambda2}.

\begin{defn} We say that a differential $1$-form $\omega$ is {\em dicritically adjusted} to $\mathcal C$ if and only if $\nu_E(\omega)=t_{s+1}$ and $\nu_{\mathcal C}(\omega)=\infty$. An {\em extended standard basis} for $\mathcal C$ is a list
$$
\omega_{-1}=dx, \omega_0=dy, \omega_{1},\ldots,\omega_{s};\, \omega_{s+1}
$$
where $\omega_{-1}, \omega_0, \omega_{1},\ldots,\omega_{s}$ is a standard basis and $\omega_{s+1}$ is
dicritically adjusted to $\mathcal C$.
\end{defn}

\begin{lema}
\label{lema:existenciaformaajustada}
 Assume that $\nu_E(\eta)=t_{s+1}$ and $\nu_\mathcal C(\eta)>u_{s+1}$. Then, there is a $1$-form $\tilde \eta$ such that $\nu_E(\tilde\eta)=t_{s+1}$ and $\nu_\mathcal C(\tilde\eta)>\nu_\mathcal C(\eta)$.
\end{lema}
\begin{proof} Take a standard basis $dx,dy, \omega_{1},\ldots,\omega_{s}$ for $\mathcal C$.   There is an index $k$ such that $\nu_\mathcal C(\eta)=an+bm+\lambda_k$. Consider the $1$-form
$
x^ay^b\omega_k
$. Note that $\nu_\mathcal C(x^ay^b\omega_k)=\nu_{\mathcal C}(\eta)$.
If we show that $an+bm+t_k>t_{s+1}$, we are done, by taking $\tilde \eta=\eta-\mu x^ay^b\omega_k$ for a convenient constant $\mu$. We have
\begin{equation}
an+bm+\lambda_k> u_{s+1}\Rightarrow an+bm+t_k> u_{s+1}-\lambda_k+t_k.
\end{equation}
It remains to show that $ u_{s+1}-\lambda_k+t_k\geq t_{s+1}$. We have
$$
u_{s+1}-\lambda_k+t_k\geq t_{s+1}\Leftrightarrow
u_{s+1}-\lambda_k+t_k\geq t_s+u_{s+1}-\lambda_s\Leftrightarrow
\lambda_s-\lambda_k\geq t_s-t_k.
$$
We are done by Lemma \ref{lema:lambdamenoslambda}.
\end{proof}
\begin{prop}
\label{prop:existenciaformaajustada}
 There is at least one $1$-form $\omega$ dicritically adjusted to $\mathcal C$.
\end{prop}
\begin{proof} Take a standard basis $dx,dy, \omega_{1},\ldots,\omega_{s}$ for $\mathcal C$.  There is an index $k<s$ such that
$
u_{s+1}=an+bm+\lambda_s=cn+dm+\lambda_k
$.
Note that
$$
an+bm+t_s<cn+dm+t_k,
$$
since $t_s-t_k<\lambda_s-\lambda_k$. In this way, we have
\begin{enumerate}
\item $t_{s+1}=\nu_E(x^ay^b\omega_s)< \nu_E(x^cy^d\omega_k)$.
\item $u_{s+1}=\nu_\mathcal C(x^ay^b\omega_s)=\nu_\mathcal C(x^cy^d\omega_k)$.
\end{enumerate}
Taking $\eta=x^ay^b\omega_s-\mu x^cy^d\omega_k$, for a convenient constant $\mu$, we have that $$\nu_E(\eta)=t_{s+1},\quad \nu_\mathcal C(\eta)>u_{s+1}.$$
By a repeated application of Lemma \ref{lema:existenciaformaajustada}, we find a $1$-form $\tilde\eta$ such that
$$\nu_E(\tilde \eta)=t_{s+1}<\nu_\mathcal C(\tilde\eta),\quad \nu_\mathcal C(\tilde\eta)>c_\Gamma+1.$$
Now, we can integrate $\tilde\eta$ as follows. Let
$\phi:t\mapsto \phi(t)=(\phi_1(t),\phi_2(t))$
be a reduced parametrization of the curve $\mathcal C$. Then $\phi$ induces a morphism
$$
\phi^{\#}:\mathbb C\{x,y\}\rightarrow \mathbb C\{t\},\quad h\mapsto h\circ\phi,
$$
whose kernel is generated by a local equation of $\mathcal C$ and moreover, the conductor ideal
$
t^{c_\Gamma}\mathbb C\{t\}
$
is contained in the image of $\phi^{\#}$. Let us write
$$
\phi^*\tilde\eta=\xi(t)dt,\quad  \xi(t)\in t^{c_\Gamma}\mathbb C\{t\}.
$$
By integration, there is a series $\psi(t)$ such that $\psi'(t)=\xi(t)$, with $\psi(t)\in t^{c_\Gamma+1}\mathbb C\{t\}$. In view of the properties of the conductor ideal, there is a function $h\in \mathbb C\{x,y\}$ such that $h\circ\phi(t)=\psi(t)$. If we consider
$
\omega=\tilde\eta-dh
$,
we have that $\nu_E(\omega)=\nu_E(\tilde\eta)=t_{s+1}$ and $\nu_\mathcal C(\omega)=\infty$.
\end{proof}
\begin{prop} If $\omega$ is dicritically adjusted to $\mathcal C$, we have that $\omega$ is basic and resonant (hence it is $E$-totally dicritical) and $\mathcal C$ is an $\omega$-cusp.
\end{prop}
\begin{proof} Recall that $t_{s+1}<nm$ and $\nu_\mathcal C(\omega)=\infty>t_{s+1}$.
\end{proof}

\subsection{Delorme's decompositions}
Delorme's decompositions are described in next Theorem \ref{th:delormedecomposition}. This result is the main tool we need to use in our statements on the analytic semiroots. We provide a proof of it, using a different approach to the one of Delorme, in Appendix~\ref{apendice-Delorme-decomp}.

\begin{theorem}[Delorme's decomposition]
\label{th:delormedecomposition} Let $\mathcal C$ be a cusp in $\operatorname{Cusps}(E)$, consider
an  extended standard basis $\omega_{-1},\omega_0,\omega_1,\ldots,\omega_s;\omega_{s+1}$ of $\mathcal C$ and denote by
$$
\Lambda=\Gamma(n,m,\lambda_1,\ldots, \lambda_s)
$$
the semimodule of differential values of $\mathcal C$. For any indices $0\leq j\leq i\leq s$, there is a decomposition
$$
\omega_{i+1}=\sum_{\ell=-1}^{j}f_\ell^{ij}\omega_\ell
$$
such that, for any $-1\leq \ell\leq j $ we have $\nu_{\mathcal C}(f_\ell^{ij}\omega_\ell)\geq v_{ij}$, where
$
v_{ij}=t_{i+1}-t_j+\lambda_j
$
and there is exactly one index $-1\leq k\leq j-1$ such that
$
\nu_{\mathcal C}(f_k^{ij}\omega_k)=\nu_{\mathcal C}(f_j^{ij}\omega_j)= v_{ij}
$.
\end{theorem}
\begin{remark}
\label{rk:cotavij}
 Note that $v_{ii}= t_{i+1}-t_i+\lambda_i=u_{i+1}$. We also have that
$$
v_{ij}=u_{i+1}-(\lambda_i-u_i)-\cdots-(\lambda_{j+1}-u_{j+1}).
$$
In particular we have that $v_{ij}\leq u_{i+1}<c_\Gamma+n< nm$. See Lemma \ref{lema:cotaconductor}.
\end{remark}

\subsection{Analytic Semiroots} Let  $\Lambda=\Gamma(n,m,\lambda_1,\ldots,\lambda_s)$ be the semimodule of differential values of the $E$-cusp $\mathcal C$ and let us consider an extended standard basis
$$
\mathcal E=(\omega_{-1},\omega_0,\omega_1,\ldots,\omega_s;\omega_{s+1})
$$

We recall that $\omega_1,\omega_2,\ldots,\omega_{s+1}$ are basic and resonant. Fix a free point $P\in E$. For each $i=1,2,\ldots,s+1$, there is an  $E$-cusp $\mathcal C^i_P$ passing through $P$ and invariant by $\omega_i$. Let us note that if $P$ is the infinitely near point of $\mathcal C$ in $E$, we have that $\mathcal C^{s+1}_P=\mathcal C$.
\begin{defn} The cusps $\mathcal C^i_P$, for $i=1,2,\ldots,s+1$ are called the {\em analytic semiroots} of $\mathcal C$ through $P$ with respect to the extended standard basis $\mathcal E$.
\end{defn}

Let us denote $\mathcal E_i=(\omega_{-1},\omega_0,\omega_1,\ldots,\omega_{i-1};\omega_{i})$ for any $1\leq i\leq s+1$. The main objective of this paper is to show the following Theorem
\begin{theorem}
\label{ref:teomain}
 For any $1\leq i\leq s+1$ and any free point $P\in E$, we have
\begin{enumerate}
\item [a)] $\Lambda_{i-1}$ is the semimodule of differential values of $\mathcal C^i_P$.
\item [b)] $\mathcal E_i$ is an extended standard basis for $\mathcal C^i_P$.
\end{enumerate}
\end{theorem}

Let us consider an index $1\leq i\leq s+1$ and an analytic semiroot $\mathcal C^i_P$ in order to prove Theorem \ref{ref:teomain}.

\begin{lema}
\label{lema:tresfinal}
We have that $\nu_{{\mathcal C}^i_P}(\omega_j)=\lambda_j$, for any $j=-1,0,1,\ldots,i-1$.
\end{lema}
\begin{proof} For any basic non resonant $1$-form $\omega$, we have that
$$
\nu_\mathcal C(\omega)=\nu_E(\omega)=\nu_{{{\mathcal C}^i_P}}(\omega).
$$
This is particularly true for the case of $1$-forms of the type $\omega=df$, where we know that $\nu_E(\omega)=\nu_E(f)$. Hence, for any germ of function $f\in \mathbb C\{x,y\} $ we have that
\begin{equation}
\min\{\nu_{\mathcal C}(df),nm\}=\min\{\nu_{{{\mathcal C}^i_P}}(df),nm\}.
\end{equation}
The statement of the Lemma is true for $\ell=-1,0$, since
$$
\nu_{\mathcal C}(dx)= \nu_{{{\mathcal C}^i_P}}(dx)=n,\quad
\nu_{\mathcal C}(dy)= \nu_{{{\mathcal C}^i_P}}(dy)=m.
$$
Let us assume that it is true for any $\ell=-1,0,1,\ldots,j$, with $0\leq j\leq i-2$. We apply
 Theorem \ref{th:delormedecomposition} to obtain a decomposition
$$
\omega_{i}=\sum_{\ell=-1}^{j+1}f_\ell\omega_\ell
$$
such that $\nu_{\mathcal C}(f_\ell\omega_\ell)\geq v_{i-1,j+1}$, where $v_{i-1,j+1}<nm$ (see Remark \ref{rk:cotavij}) and there is a single $k\leq j$ such that
$
\nu_{\mathcal C}(f_{j+1}\omega_{j+1})=\nu_{\mathcal C}(f_k\omega_k)=v_{i-1,j+1}.
$
We deduce that
$$
\nu_{\mathcal C}\left(\sum_{\ell=-1}^{j}f_\ell\omega_\ell\right)=\nu_{\mathcal C}(f_k\omega_k)=v_{i-1,j+1}.
$$
On the other hand, by induction hypothesis and noting that $v_{i-1,j+1}<nm$, we have that
$$
\min\{\nu_{{\mathcal C}}(f_\ell\omega_\ell),nm\}=\min\{\nu_{{{\mathcal C}^i_P}}(f_\ell\omega_\ell),nm\},\quad \ell=-1,0,1,\ldots,j.
$$
In particular, we have that
$$\nu_{{{\mathcal C}^i_P}}(f_k\omega_k)=v_{i-1,j+1},\quad
\nu_{{{\mathcal C}^i_P}}\left(\sum_{\ell=-1}^{j}f_\ell\omega_\ell\right)=v_{i-1,j+1}.
$$
Since $\nu_{\mathcal C^i_P}(\omega_{i})=\infty$, we have that
$
\nu_{{{\mathcal C}^i_P}}(f_{j+1}\omega_{j+1})=v_{i-1,j+1}$. Hence we have
$$
v_{i-1,j+1}=\nu_{{\mathcal C}}(f_{j+1}\omega_{j+1})=\nu_{{{\mathcal C}^i_P}}(f_{j+1}\omega_{j+1}).
$$
Noting that $\nu_{\mathcal C}(f_{j+1})=\nu_{{{\mathcal C}^i_P}}(f_{j+1})$, we conclude that
$
\nu_{\mathcal C}(\omega_{j+1})=\nu_{{{\mathcal C}^i_P}}(\omega_{j+1})
$.
The proof is ended.
\end{proof}
\begin{corolario} $\Lambda^{{{\mathcal C}^i_P}}\supset\Lambda_{i-1}$.
\end{corolario}
\begin{proof} It is enough to remark that $\lambda_j\in \Lambda^{{{\mathcal C}^i_P}}$ for any $j=-1,0,1,\ldots,i-1$.
\end{proof}

\begin{prop}
\label{pro:igualdaddesemimodulos2} $\Lambda^{{\mathcal C}^i_P}=\Lambda_{i-1}$.
\end{prop}
\begin{proof} We already know that $\Lambda^{{{\mathcal C}^i_P}}\supset\Lambda_{i-1}$. Assume that $\Lambda^{{{\mathcal C}^i_P}}\ne\Lambda_{i-1}$ and take the  number
$$
\lambda=\min\left(\Lambda^{{{\mathcal C}^i_P}}\setminus \Lambda_{i-1}\right).
$$
Note that $\lambda>m$ and hence there is a maximum index $0\leq j\leq i-1$ such that
$
\lambda_j<\lambda
$. We have that
$$
\nu_E(\omega_j)=t_j=t_j(\mathcal C^i_P),\quad   t_j=t_j(\mathcal C).
$$
where $t_j(*)$ denotes the critical divisorial order $t_j$ with respect to the curve $*$. Assume first that $j\leq i-2$.

Let $\tilde\omega$ be a $1$-form in a standard basis for ${{\mathcal C}^i_P}$ that corresponds to the differential value $\lambda$. Note that $\lambda$ is the differential value in the basis of $\Lambda^{\mathcal C^i_P}$ that immediately follows $\lambda_j$ and the previous ones correspond to the values in the basis of $\Lambda$. We have that
$$
\nu_E(\omega_{j+1})=t_{j+1}=t_j+u_{j+1}-\lambda_j=t_{j+1}(\mathcal C^i_P)=\nu_E(\tilde\omega).
$$
In view of the property (1) in Theorem \ref{teo:bs:delorme}, we have that
$$
\lambda=\max\{\nu_{{{\mathcal C}^i_P}}(\omega);\; \nu_E(\omega)=t_{j+1}(\mathcal C^i_P)=t_{j+1}\}
$$
In view of Lemma \ref{lema:tresfinal},
we know that $\nu_{{{\mathcal C}^i_P}}(\omega_{j+1})=\lambda_{j+1}$ and $\nu_E(\omega_{j+1})=t_{j+1}$.
This implies that $\lambda_{j+1}<\lambda$, contradiction.

Let us consider now the case when $j=i-1$. We shall prove that it is not possible to have $s(\mathcal C^i_P)>i-1$ where $s(*)$ refers to ``concept $s$'' with respect to the curve $*$ (that is, $s+2$ is the number of elements of the basis of the semimodule of the curve). If $s(\mathcal C^i_P)>i-1$, we have that
$$
\lambda=\max\{\nu_{{{\mathcal C}^i_P}}(\omega);\; \nu_E(\omega)= t_{i}(\mathcal C^i_P)=t_{i}\}.
$$
But we know that $\nu_{{{\mathcal C}^i_P}}(\omega_{i})=\infty$ and $\nu_E(\omega_{i})=t_{i}$, this is the desired contradiction.
\end{proof}

\begin{example}
Let us compute the analytic semiroots of the curve $\mathcal{C}$ given in Example~\ref{ex:5-11-v2}. The Puiseux parametrizations  for the $E$-cusps $\mathcal{C}_a^i$ of $\omega_i$, $i=1,2,3$, are given by
\begin{align*}
\phi_1^a(t)& =(t^5,at^{11}) \\
\phi_2^a(t) & = (t^5, a t^{11}+ a^2 t^{12} + \tfrac{23}{22}a^3 t^{13}+ \tfrac{136}{121} a^4 t^{14} + \cdots) \\
\phi_3^a(t) & = (t^5, \sum_{i \geq 11} a^{i-10} t^{i})
\end{align*}
with $a \in \mathbb{C}^*$. Hence, the analytic semiroots of $\mathcal{C}$ are the curves $\mathcal{C}^1_1$, $\mathcal{C}^2_1$ and $\mathcal{C}^3_1$ given by the above parametrizations $\phi_i^1(t)$, $i=1,2,3$, with $a=1$.

Note that in this example, for any  $i=1,2,3$,  all the $E$-cusps of the family $\{\mathcal{C}^i_a\}_a$  are analytically equivalent. To see this it is enough to consider the new parameter $t=a^{-1}u$ and the change of variables $x_1=a^5x$, $y_1=a^{10} y$.
\end{example}
\begin{example}\label{ex:7-17}
We would like to remark that, in general, the $E$-cusps of an element $\omega_i$, $i \geq 2$, of a standard basis are not analytically equivalent as the following example shows. Consider the curve $\mathcal{C}$ given by the Puiseux parametrization \linebreak $\phi(t)=(t^7,t^{17}+t^{30}+t^{33}+t^{36})$ with $\Gamma=\langle 7,17\rangle$ and $\Lambda=\Gamma(7,17, 37,57)$. A standard basis is given by
$\omega_{-1}=dx, \quad \omega_0=dy, \quad \omega_1=7xdy-17 y dx$
and
$$\omega_2=3757x^2ydx-1547x^3dy-4624y^2dx+1904xydy+1183 y^2dy.$$
The $E$-cusps of $\omega_2$ are the curves given by the Puiseux parametrization
$$\varphi_a(t)=(t^7, at^{17} + a^3 t^{30} + a^4t^{33}+ \cdots)$$
with $a \in \mathbb{C}^*$. If we consider a new parameter $t=a^{-2/13}u$ and the we make the change of variables $x_1=a^{14/13}x, y_2=a^{21/13} y$, we obtain that the family of $E$-cusps of $\omega_2$ are the curves $\mathcal{C}^2_a$ given by the parametrizations
$$\phi_a(t)=(t^7,t^{17}+t^{30}+a^{7/13} t^{33} + \cdots)$$
From the results above, we have that $\Lambda^{\mathcal{C}_a^2}=\Lambda_1=\Gamma(7,17,37)$ for all $a \in \mathbb{C}^*$. Since $33 \not \in \Lambda_1 - 7$, by Theorem 2.1  in \cite{hefez2}, two curves  $\mathcal{C}^2_{a_1}$ and $\mathcal{C}^2_{a_2}$  are not, in general, analytically equivalent  for  $a_1,a_2 \in \mathbb{C}^*$.
\end{example}
\begin{example}
Let us consider the 1-form $\omega$ of example 4.7 in \cite{Oziel} given by
$$
\omega=(7y^5+2x^9y-2x^9y^2-9x^2y^4) dx + (4y^3x^3-x^{10} + 2x^{10} y -3xy^4-x^{8}y^2)dy.
$$
This 1-form is pre-basic and resonant for the pair $(4,9)$ since $\nu_E(\omega)=48$, the co-pair of $(4,9)$ is $(3,7)$ and $\text{Cl}(\omega;x,y) \subset R^{4,9}(3,4)$ where $$R^{4,9}(3,4)=\{(\alpha,\beta)\in \mathbb{R}^2 \ : \ \alpha + 2\beta \geq 11\} \cap \{(\alpha,\beta)\in \mathbb{R}^2 \ : \ 3 \alpha + 7 \beta \geq 37\}.$$ Moreover, the weighted initial part of $\omega$ is given by $\text{In}_{4,9;x,y}^{48}=x^2y^3 (-9ydx+4xdy)$.
Consequently $\omega$ is totally $E$-dicritical for the last divisor $E$ associated to the cuspidal sequence $\mathcal{S}^{4,9}_{y=0}$. Note that $\omega$ is not a basic 1-form since $\nu_E(\omega) > nm=36$.

The invariant curves of $\omega$ which are transversal to the dicritical component $E$ are the curves $\mathcal{C}_a$, $a \in \mathbb{C}^*$, given  by
$$y^4-ax^9+(a-1)x^7y+x^7y^2=0.$$
Note that these curves  do not have the same semimodule of differential values since  the curves $\mathcal{C}_a$, with $a \neq 1$, have   Zariski invariant equal to 10 whereas  the curve $\mathcal{C}_1$ has   Zariski invariant   equal to 19.
\end{example}
\appendix
\section{Bounds for the Conductor}

In this Appendix we present some bounds for the conductor of the semimodules $\Lambda_i$ in the decomposition series of
 $$
 \Lambda=\Gamma(\lambda_{-1},\lambda_0,\lambda_1,\ldots,\lambda_s)
 $$
 that will be enough to prove the results we need relative to the structure of the semimodule of differential values. We recall that the semigroup $\Gamma$ is generated by a Puiseux pair $(n,m)$, with $2\leq n$.

Given two integer numbers $r\leq s$, we denote $[r,s]$ the set of the integer numbers $\ell$ such that $r\leq \ell\leq s$. For any $q\geq 0$ we denote by $I_q$ the interval
$
I_q=[nq,n(q+1)-1]
$, in particular, we have that $I_0=\{0,1,2,\ldots,n-1\}$.
For any $r,s\in I_0$, we define the {\em circular interval} $<r,s>$ by
$$
<r,s>=[r,s],\; \text{ if } r\leq s;\quad <r,s>=[r,n-1]\cup[0,s],\; \text{ if } r> s.
$$
We denote by $\rho:\mathbb Z\mapsto \mathbb Z/(n)$ the canonical map and we also use the notation $\rho(p)=\bar{p}$.
Since $\gcd(n,m)=1$, there is a bijection
$$
\xi:\mathbb Z/(n)\rightarrow I_0=\{0,1,2,\ldots,n-1\},\quad  \xi^{-1}(k)= \rho(km).
$$

For any $q\geq 0$ and any subset $S \subset \mathbb{Z}_{\geq 0}$, we define the {\em $q$-level set $R_q(S)\subset I_0$} by
$
R_q(S)=\xi\left(\rho( S \cap I_q)\right)
$.

\begin{remark}\label{re:insemi:contencion}
We have that $R_q(\Lambda)\subset R_{q'} (\Lambda)$ for all $q'\geq q$.
\end{remark}

\begin{remark}\label{re:aumentar}
For any $\mu\in \mathbb Z_{\geq 0}$ and $q\geq 1$, we have
$$
\#R_{q}(\mu+\Gamma)\leq\#R_{q-1}(\mu+\Gamma)+1.
$$
Indeed, this is equivalent to show that $\#\rho((\mu+\Gamma)\cap I_q)\leq\#\rho((\mu+\Gamma)\cap I_{q-1})+1$. Assume that $\bar p_1,\bar p_2\in \rho((\mu+\Gamma)\cap I_q)\setminus \rho((\mu+\Gamma)\cap I_{q-1})$. We can take representatives $p_1,p_2\in (\mu+\Gamma)\cap I_q$ of $\bar p_1$ and $\bar p_2$ of the form
$
p_1=\mu+b_1m$, $p_2=\mu+b_2m
$.
If $p_1\ne p_2$, we have that $|p_1-p_2|\geq m>n$ and this is not possible.
\end{remark}
\begin{lema}\label{le:insemi:circular} Consider $\mu\in I_v$, denote $r=\xi(\bar\mu)$ and let $q$ be such that $q\geq v$.  For any $p\in R_q(\mu+\Gamma)$ we have that $<r,p>\subset R_q(\mu+\Gamma)$. In particular, the set
 $R_q(\mu+\Gamma)$ is a circular interval.
\end{lema}
\begin{proof} The second statement is straightforward, since the union of circular intervals with a common point is a circular interval. To prove the first statement, we
proceed by induction on the number $\ell$ of elements in $<r,p>$. If $\ell\leq 2$, we are done, since $<r,p>\subset\{r,p\}\subset R_q(\mu+\Lambda)$. Assume that $\ell>2$; in particular we have that $r\ne p$.  Consider the point $\tilde p\in I_0$ given by $\tilde p=p-1$, if $p\geq 1$ and $\tilde p=n-1$, if $p=0$.
We have that $<r,p>=<r,\tilde p>\cup\{p\}$ and the length of $<r,\tilde p>$ is $\ell-1$. Then, it is enough to show that $\tilde p\in R_q(\mu+\Gamma)$.  Take an element
$
\mu+an+bm\in I_q\cap (\mu+\Gamma)
$
such that $\rho({\mu+an+bm})=\rho({pm})$.  Noting that $r\ne p$, we have that $b\geq 1$. There is $q'\leq q$ such that $\mu+an+(b-1)m\in I_{q'}$ and hence
$$
\mu+(a+q-q')n+(b-1)m\in I_{q}\cap (\mu+\Gamma).
$$
We have that
$
\rho({\mu+(a+q-q')n+(b-1)m})=\rho({\tilde p m})
$
and thus $\tilde p\in R_q(\mu+\Gamma)$.
\end{proof}
\begin{defn}
We define the {\em tops $q_1$ and $q_2$ of $\Lambda$} by the property that
$$
\lambda_s+n\ell_1\in I_{q_1}, \quad \lambda_s+m\ell_2\in I_{q_2}
$$
where $\ell_1$ and $\ell_2$ are the limits of $\Lambda$. The {\em main top $Q_\Lambda$} is the maximum of $q_1,q_2$.
\end{defn}
\begin{prop}\label{lema:conductor1}
Consider a normalized semimodule $\Lambda=\Gamma(0,\lambda_0,\lambda_1,\ldots,\lambda_s)$. Let $v$ be such that $\lambda_s\in I_v$ and assume that $R_q(\Lambda_{s-1})$ is a circular interval for any $q\geq v$. Denote by  $q_1,q_2$ the tops of $\Lambda$ and put $r=\xi(\bar\lambda_s)$.
Then we have that
\begin{enumerate}
 \item $[0,r-1]\subset R_q(\Lambda_{s-1})$,   for all $q\geq q_1-1$.
 \item $[r,n-1]\subset R_q(\Lambda)$,  for all $q\geq q_2-1$.
\end{enumerate}
In particular $c_\Lambda\leq n(Q_\Lambda-1)$.
\end{prop}
\begin{proof} Note that Statements (1) and (2) imply that
$$
I_0=[0,n-1]=[0,r]\cup [r,n-1]\subset R_q(\Lambda),\quad q\geq Q_\Lambda-1
$$
and thus, we have that $c_\Lambda\leq n(Q_\Lambda-1)$.

{\em Proof of Statement (1):} By Remark \ref{re:insemi:contencion} it is enough to show that  we have $[0,r-1]\subset R_{q_1-1}(\Lambda_{s-1})$. Since $\lambda_s+n\ell_1\in \Lambda_{s-1}$, there is an index $k\leq s-1$ such that
$
\lambda_s+n\ell_1=\lambda_k+an+bm
$. By the minimality of $\ell_1$, we have that $a=0$ and hence
$
\lambda_s+n\ell_1=\lambda_k+bm
$.
Denote $r_k=\xi(\bar\lambda_k)$. Note that $r_k\ne r$, since $r\notin R_v(\Lambda_{s-1})$.

Assume that the next statements are true:
\begin{enumerate}
\item[a)] If $r_k>r$, then $[0,r]\subset R_{q_1}(\lambda_k+\Gamma)$.
\item[b)] If $r_k<r$, then $[r_k,r]\subset R_{q_1}(\lambda_k+\Gamma)$ and $[0,r_k]\subset R_{q_1-1}(\Lambda_{s-1})$.
\end{enumerate}
If $r_k>r$, by the minimality of $\ell_1$, we have that $r\notin R_{q_1-1}(\lambda_k+\Gamma)$; now,  in view of Remark \ref{re:aumentar} and noting that $[0,r]=[0,r-1]\cup \{r\}$,  we obtain that
$$[0,r-1]\subset R_{q_1-1}(\lambda_k+\Gamma)\subset R_{q_1-1}(\Lambda_{s-1}).$$
If $r_k<r$, we obtain as above that $[r_k,r-1]\subset  R_{q_1-1}(\lambda_k+\Gamma)$, then
$$
[0,r-1]=[0,r_k]\cup [r_k,r-1]\subset R_{q_1-1}(\Lambda_{s-1}).
$$
If remains to prove a) and b).

{\em Proof of a):} We can apply Lemma \ref{le:insemi:circular} to have that $<r_k,r>\subset R_{q_1}(\lambda_k+\Gamma)$. We end by noting that $[0,r]\subset <r_k,r>$.

{\em Proof of b):} We apply Lemma \ref{le:insemi:circular} to have that $<r_k,r>=[r_k,r]\subset R_{q_1}(\lambda_k+\Gamma)$. On the other hand, we know that $R_{q_1-1}(\Lambda_{s-1})$ is a circular interval since $q_1-1\geq v$ and it contains $0$ and $r_k$. Moreover $r\notin R_{q_1-1}(\Lambda_{s-1})$ and $r>r_k$, then the circular interval $R_{q_1-1}(\Lambda_{s-1})$ contains $[0,r_k]$.

{\em Proof  of Statement (2):}
 It is enough to show that $[r,n-1]\subset R_{q_2-1}(\Lambda)$. By an argument as before, there is an index $k\leq s-1$ such that
$
\lambda_s+m\ell_2=\lambda_k+na
$. Take $r_k\ne r$ as above. By Lemma \ref{le:insemi:circular}, we have that
$
<r,r_k>\subset R_{q_2}(\lambda_s+\Gamma)
$.
Let us see that $r_k\notin R_{q_2-1}(\lambda_s+\Gamma)$. For this, let us show that the  property
$$r_k\in R_{q_2-1}(\lambda_s+\Gamma)$$ leads to a contradiction. This property should imply that $\lambda_k+n(a-1)\in \lambda_s+\Gamma$ and hence there are nonnegative integer numbers $\alpha,\beta$ such that
$$
\lambda_s+n\alpha+m\beta=
\lambda_k+n(a-1).
$$
If $a-1\leq \alpha$, we have that $\lambda_k=\lambda_s+(\alpha-a+1)n+\beta m$ and this contradicts the fact that $\lambda_k<\lambda_s$; hence  $a-1>\alpha$ and we have
$$
\lambda_s+m\beta=
\lambda_k+n(a-1-\alpha).
$$
Since $a-1-\alpha<a$, we have that $\beta<\ell_2$. This contradicts the minimality of $\ell_2$.

Since $r_k\notin R_{q_2-1}(\lambda_s+\Gamma)$, we can apply Remark \ref{re:aumentar} that tells us that
$$
<r,r_k>\setminus \{r_k\}\subset R_{q_2-1}(\lambda_s+\Gamma)\subset R_{q_2-1}(\Lambda).
$$
Note also that $r_k\in R_{q_2-1}(\Lambda)$. Then we have that
$
<r,r_k>\subset R_{q_2-1}(\Lambda)
$.

If $r>r_k$, then $[r,n-1]\subset <r,r_k>\subset R_{q_2-1}(\Lambda) $. Assume now that $r<r_k$. Recall that $r\notin R_v(\Lambda_{s-1})$; since $R_v(\Lambda_{s-1})$ is a circular interval containing $r_k$ and $0$, but not containing $r$, we have that
$$[r_k,n-1]\subset R_v(\Lambda_{s-1})\subset R_{q_2-1}(\Lambda_{s-1})\subset R_{q_2-1}(\Lambda).
$$
We conclude that
$
[r,n-1]=<r,r_k>\cup [r_k,n-1]\subset  R_{q_2-1}(\Lambda)
$.
\end{proof}

\begin{prop}\label{prop:insemi:creciente}
Let $\Lambda$ be a normalized  increasing semimodule of length $s$ and let $v$ be such that $u_{s+1}\in I_v$. Then $R_q(\Lambda)$ is a circular interval for any $q\geq v$.
\end{prop}

\begin{proof}
Let us proceed by induction on the length $s$ of $\Lambda$. If $s=-1$, we have $\Lambda=\Lambda_{-1}=\Gamma$. By Lemma
\ref{le:insemi:circular} applied to $\mu=0$, we are done. Let us suppose that $s\geq 0$ and assume by induction that the result is true for $\Lambda_{s-1}$.
We have that $\Lambda=\Lambda_{s-1}\cup (\lambda_s+\Gamma)$. This implies that
$$
R_q(\Lambda)=R_q(\Lambda_{s-1})\cup R_q(\lambda_s+\Gamma),\quad q\geq0.
$$
Let $v'$ be such that $u_{s}\in I_{v'}$. By induction hypothesis, we know that
$R_q(\Lambda_{s-1})$ is a circular interval for any $q\geq v'$. Moreover, by the increasing property, we have that
$$
u_{s+1}>\lambda_s>u_s\geq \lambda_{s-1}.
$$
In particular, we have that $v\geq v'$ and
$R_q(\Lambda_{s-1})$ is a circular interval for any $q\geq v$.
On the other hand, take $v''$ such that $\lambda_s\in I_{v''}$. By Lemma \ref{le:insemi:circular}, we know  that $R_q(\lambda_s+\Gamma)$ is a circular interval for any $q\geq v''$. Since $v\geq v''$, we have that $R_q(\lambda_s+\Gamma)$ is a circular interval for any $q\geq v$. Thus, both
$$
R_q(\Lambda_{s-1})\mbox{ and } R_q(\lambda_s+\Gamma)
$$
are circular intervals for $q\geq v$. We need to show that their union is also a circular interval.  Since $v\geq v''$, we have that $r\in R_q(\lambda_s+\Gamma)$ for
$r=\xi(\overline{\lambda}_s)$. Noting that $0\in R_q(\Lambda_{s-1})$ and $r\in R_q(\lambda_s+\Gamma)$, in order to prove that $R_q(\Lambda)$ is a circular interval, it is enough to show that one of the following properties holds
\begin{equation*}
a):\quad  [0,r]\subset R_q(\Lambda); \quad\quad
b):\quad  [r,n-1]\subset R_q(\Lambda).
\end{equation*}
We can apply Proposition \ref{lema:conductor1} to $\Lambda=\Lambda_{s-1}(\lambda_s)$. Indeed, by induction hypothesis we know that $R_{p}(\Lambda_{s-1})$ is a circular interval for any $p\geq v'$; since $v'\leq v''$, we have that $R_{p}(\Lambda_{s-1})$ is a circular interval for any $p\geq v''$ and thus we are in the hypothesis of Proposition \ref{lema:conductor1}. Now, by Equation \eqref{eq:uesemasuno} in Remark \ref{rk:expresionesui}, we have either $u_{s+1 }=\lambda_s+nl_1$ or $u_{s+1}=\lambda_s+ml_2$.

a)  If we have $u_{s+1 }=\lambda_s+nl_1$, noting that $q_1=v$,  we apply
Proposition \ref{lema:conductor1} and we obtain that
$
[0,r]\subset R_q(\Lambda_{s-1})\subset R_q(\Lambda)
$.

b)  If we have $u_{s+1 }=\lambda_s+ml_2$, noting that $q_2=v$,  we apply
Proposition \ref{lema:conductor1} and we obtain that $
[r,n-1]\subset  R_q(\Lambda)
$.
\end{proof}
\begin{corolario}
 \label{cor:conductor}
 Let $\Lambda$ be an increasing semimodule such that its minimum element $\lambda_{-1}$ is a multiple of $n$. Let $Q_\Lambda$ be the main top of $\Lambda$. Then
$
  c_\Lambda\leq n(Q_\Lambda-1)
$.
\end{corolario}
\begin{proof} Assume first that $\Lambda$ is normalized. Let $v'$ be such that $u_{s}\in I_{v'}$ and $v''$ such that $\lambda_s\in I_{v''}$. We know that $v''\geq v'$.  By Proposition  \ref{prop:insemi:creciente}, we know that $R_q(\Lambda_{s-1})$ is a circular interval for any $q\geq v'$ and hence for any $q\geq v''$. Then we are in the hypothesis of Proposition \ref{lema:conductor1} and we conclude.

Assume now that $\lambda_{-1}=kn$ and consider the normalization $\widetilde{\Lambda}=\Lambda-kn$. Let us note that the tops are related by the property
$
\tilde q_j=q_j-k$, for $j=1,2$
 and hence $Q_{\widetilde{\Lambda}}=Q_{\Lambda}-k$. On the other hand $c_{\widetilde{\Lambda}}=c_{\Lambda}-nk$.
 We conclude that
 $$c_{\Lambda}=c_{\widetilde{\Lambda}}+nk\leq n(Q_{\widetilde{\Lambda}}-1)+nk=
 n(Q_{\Lambda}-1).$$
\end{proof}

\section{Structure of the Semimodule}  \label{Apendice-Semimodulo}
 In this Appendix  we present a proof, using a different approach to the one of Delorme, of the main results on the structure of the semimodule of differential values for an $E$-cusp $\mathcal C$.  As before, we denote
 $
 \Lambda=\Gamma(n,m,\lambda_1,\ldots,\lambda_s)$, $ n\geq 2
 $,
the semimodule of differential values of $\mathcal C$ and we select a standard basis
$
\omega_{-1}=dx,\omega_0=dy,\omega_1,\ldots,\omega_s
$ of the cusp $\mathcal C$.

\begin{prop}\label{prop:bs:delorme}
For each $1\leq i\leq s$ we have the following statements
\begin{enumerate}
\item $\lambda_i=\sup\{\nu_\mathcal C(\omega):\nu_{E}(\omega)=t_i\}$.
\item If $\nu_\mathcal C(\omega)=\lambda_i$, then $\nu_{E}(\omega)=t_i$.
\item For each $1$-form $\omega$ with $\nu_\mathcal C(\omega)\notin \Lambda_{i-1}$, there is a unique pair $a,b\geq 0$ such that $\nu_{E}(\omega)=\nu_{E}(x^ay^b\omega_i)$. Moreover, we have that $\nu_\mathcal C(\omega)\geq \lambda_i+na+mb$.
\item We have that $\lambda_i>u_i$.
\item Let $k=\lambda_i+na+mb$, then $k\notin \Lambda_{i-1}$ if and only if for all $\omega$ such that $\nu_\mathcal C(\omega)=k$ we have that $\nu_{E}(\omega)\leq \nu_{E}(x^ay^b\omega_i)$.
\end{enumerate}
In particular, the semimodules $\Lambda_i$ are increasing, for $i=1,2,\ldots,s$.
\end{prop}
\begin{proof} Assume that $i=1$ and then $t_1=n+m=u_1$. We have
\begin{itemize}
\item Statement (1) is proven in Proposition \ref{prop:ordmonzariski}.
\item Statement (2) is proven in Corollary \ref{cor:ordlambdauno}.
\item Statement (3) is proven in Lemma \ref{le:bs:w15}.
\item Statement (4) follows from the fact that $\lambda_1>n+m=\nu_E(\omega_1)=u_1$.
\item Statement (5) is proven in Lemma \ref{le:bs:w14}.
\end{itemize}
Now, let us assume that $i\geq 2$ and take the induction hypothesis that the statements (1)-(5) are true for indices $\ell$ with  $1\leq\ell\leq i-1$.

Denote by $\ell_1$ and $\ell_2$ the $\Lambda_{i-1}$-limits. By Equation \eqref{eq:uesemasuno} in Remark \ref{rk:expresionesui}, we have two possibilities: either $u_i=\lambda_{i-1}+n\ell_1$ or $u_i=\lambda_{i-1}+m\ell_{2}$. We assume that $u_i=\lambda_{i-1}+n\ell_1$, the computations in the case $u_i=\lambda_{i-1}+m\ell_2$ are similar ones.

The proof is founded in three claims as follows:

$\bullet$
{\em Claim 1: There is a $1$-form $\eta$ with $\nu_E(\eta)=t_i$, whose initial part is proportional to the initial part of $x^{\ell_1}\omega_{i-1}$ and such that either $\nu_\mathcal C(\eta)\geq c_\Gamma$ or $\nu_\mathcal C(\eta)\notin \Lambda_{i-1}$. }
\vspace{5pt}

$\bullet$ {\em Claim 2: Any $1$-form $\omega$ with $\nu_{\mathcal C}(\omega)\notin \Lambda_{i-1}$ is reachable from $x^{\ell_1}\omega_{i-1}$.}
\vspace{5pt}

$\bullet$ {\em Claim 3:  Let $\eta$ be a $1$-form such that $\nu_E(\eta)=t_i$ whose initial part is proportional to the initial part of $x^{\ell_1}\omega_{i-1}$ and such that either $\nu_\mathcal C(\eta)\geq c_\Gamma$ or $\nu_\mathcal C(\eta)\notin \Lambda_{i-1}$. Then $\nu_\mathcal C (\eta)=\lambda_i$}.
\vspace{5pt}

\noindent We recall to the reader that the notion  ``initial part'' refers to the concept of weighted initial part defined in section~\ref{sec:weighted-initial-part}.

\medskip
{\em Proof of Claim 1:}
Recall that
$
t_i=\nu_E(\omega_{i-1})+u_i-\lambda_{i-1}=\nu_E(\omega_{i-1})+n\ell_1
$. Let us start with $\eta_1=x^{\ell_1}\omega_{i-1}$. We have that
$$\nu_E(\eta_1)=n\ell_1+\nu_E(\omega_{i-1})=t_i,\quad \nu_{\mathcal C}(\eta_1)=n\ell_1+\lambda_{i-1}=u_i\in \Lambda_{i-2}.
$$
By Statement (5) applied to $\nu_\mathcal C(\eta_1)\in \Lambda_{i-2}$, there is $\eta_1'$ with $\nu_\mathcal C(\eta'_1)=\nu_{\mathcal C}(\eta_1)$ and
$\nu_E(\eta'_1)>\nu_E(\eta_1)$. Since $\nu_{\mathcal C}(\eta'_1)=\nu_{\mathcal C}(\eta_1)$, there is a non-null constant $\mu$ such that
$$
\nu_{\mathcal C}(\tilde \eta)>\nu_{\mathcal C}(\eta_1),\quad \text{ where }
\tilde \eta=\eta_1-\mu\eta'_1.
$$
Since $\nu_E(\eta'_1)>\nu_E(\eta_1)$, we have that $\nu_E(\tilde\eta)=\nu_E(\eta_1)=t_i$ and the initial part of $\tilde\eta$ is the same one as the initial part of $\eta_1=x^{\ell_1}\omega_{i-1}$. If $\nu_{\mathcal C}(\tilde\eta)\geq c_\Gamma$ or
$\nu_{\mathcal C}(\tilde\eta)\notin \Lambda_{i-1}$, we put $\eta=\tilde\eta$ and we are done. Assume that $\nu_\mathcal C(\tilde\eta)\in \Lambda_{i-1}$. Let us write
$$
\nu_{\mathcal C}(\tilde\eta)=an+bm+\lambda_{\ell},\quad \ell\leq i-1.
$$
Let us see that $\nu_E(\tilde\eta)<\nu_E(x^ay^b\omega_\ell)$; this is equivalent to verify that
$
t_i-t_\ell<na+mb
$. Since $\nu_\mathcal C(\tilde\eta)>u_i$, in view of  Lemma \ref{lema:lambdamenoslambda} we have
$$
na+mb>u_i-\lambda_\ell=n\ell_1+\lambda_{i-1}-\lambda_\ell\geq n\ell_1+t_{i-1}-t_\ell=t_{i}-t_\ell.
$$
On the other hand, we have that $\nu_{\mathcal C}(\tilde\eta)=\nu_{\mathcal C}(x^ay^b\omega_{\ell})$. Thus, there is a constant $\mu_1$, such that
$\nu_\mathcal C(\tilde\eta_1)>\nu_{\mathcal C}(\tilde\eta)$, $\nu_E(\tilde\eta_1)=\nu_E(\tilde \eta)$,  where  $\tilde\eta_1=\tilde\eta-\mu_1 x^ay^b\omega_{\ell},
$
and the initial part of $\tilde\eta_1$ is the same one as the initial part of $x^{\ell_1}\omega_{i-1}$.

 If $\nu_{\mathcal C}(\tilde\eta_1)\in \Lambda_{i-1}$, we repeat the procedure starting with $\tilde\eta_1$, to obtain $\tilde\eta_2$ such that $\nu_{\mathcal C}(\tilde\eta_2)>\nu_{\mathcal C}(\tilde\eta_1)$ and $\nu_E(\tilde\eta_2)=t_i$. After finitely many repetitions, we get a $1$-form $\eta$ such that $\nu_E(\eta)=t_i$, whose initial part is the same one as $x^{\ell_1}\omega_{i-1}$ and either $\nu_{\mathcal C}(\eta)\geq c_\Gamma$ or $\nu_{\mathcal C}(\eta)\notin \Lambda_{i-1}$. This proves Claim 1.
\vspace{5pt}

{\em Proof of Claim 2: }
  Take $\omega$ such that $\lambda=\nu_{\mathcal C}(\omega)\notin \Lambda_{i-1}$. Note that $\lambda\notin \Lambda_{i-2}$.  By Statement (3), we have that $\omega$ is reachable from $\omega_{i-1}$. Thus, there are $a,b\geq 0$ and a constant $\mu$ such that
 $$
 \nu_E(\omega-\mu x^ay^b\omega_{i-1})>\nu_E(\omega)=\nu_E(x^ay^b\omega_{i-1})=an+bm+t_{i-1}.
 $$
 and moreover, we have that
 $
 \lambda=\nu_{\mathcal C}(\omega)>an+bm+\lambda_{i-1}=k
 $ (note that $\lambda\ne k$, since $\lambda\notin\Lambda_{i-1}$).

 Consider the $1$-form $\omega'=\omega-\mu x^ay^b\omega_{i-1}$. We know that
 $$
 \nu_{\mathcal C}(\omega')=k,\quad \nu_E(\omega')>\nu_E(x^ay^b\omega_{i-1}).
 $$
  By Statement (5), we conclude that $k\in \Lambda_{i-2}$. Hence $k\in \Lambda_{i-2}\cap (\lambda_{i-1}+\Gamma)$. Let us show that we necessarily have that $a\geq \ell_1$. Write
   $$
   k=an+bm+\lambda_{i-1}=\tilde a n+\tilde b m+\lambda_j,\quad j\leq i-2.
   $$
   Since $\lambda_{i-1}>\lambda_j$, we have that $an+bm<\tilde a n+\tilde b m$. Thus, we have either $a<\tilde a$ or $b<\tilde b$. If $b<\tilde b$, we have that
   $
   an+\lambda_{i-1}=\tilde a n+(\tilde b-b)m +\lambda_{j}\in \Lambda_{i-2}\cap (\lambda_{i-1}+\Gamma)
   $.
   In view of the minimality of $\ell_1$ we should have that $\ell_1\leq a$ and then  $\omega$ is reachable from
 $x^{\ell_1}\omega_{i-1}$. Assume that $a<\tilde a$ and let us obtain a contradiction. We have
   $$
   bm+\lambda_{i-1}=(\tilde a-a) n+\tilde bm +\lambda_{j}\in \Lambda_{i-2}\cap (\lambda_{i-1}+\Gamma).
   $$
   We deduce that $b\geq \ell_2$.
   By Statement (4),
    we know that $\Lambda_{i-1}$ is an increasing semimodule, starting at $\lambda_{-1}=n$.
    By Corollary  \ref{cor:conductor}, we know that
    $
   c_{\Lambda_{i-1}}\leq n(Q_{\Lambda_{i-1}}-1)
    $, where $Q_{\Lambda_{i-1}}=\max\{q_1,q_2\}$ and  $q_1,q_2$ are the tops of $\Lambda_{i-1}$.
    Suppose that $\lambda\in I_{d}$, we have
    \begin{enumerate}
    \item
    $\lambda >k=an+bm+\lambda_{i-1}\geq u_i=\ell_1n+\lambda_{i-1}$ and hence $d\geq q_1$.
    \item $\lambda >k=an+bm+\lambda_{i-1}\geq \ell_2m+\lambda_{i-1}$ and hence $d\geq q_2$.
    \end{enumerate}
    We conclude that $\lambda\in \Lambda_{i-1}$, contradiction. This ends the proof of Claim 2.
\vspace{5pt}

{\em Proof of Claim 3: }
Note that $\nu_\mathcal C(\eta)\geq\lambda_i$. Assume that $\lambda=\nu_{\mathcal C}(\eta)>\lambda_i$. Recalling that $\nu_{\mathcal C}(\omega_i)=\lambda_i\notin \Lambda_{i-1}$ and that  the initial part of $\eta$ is proportional to the initial part of $x^{\ell_1}\omega_{i-1}$, we can apply Claim 2 and we get that  $\omega_i$ is reachable from $\eta$. Then there are $a,b\geq 0$ and a constant $\mu$ such that $\nu_E(\omega_i-\mu x^ay^b\eta)>\nu_E(\omega_i)$. Put $\omega_i^1=\omega_i-\mu x^ay^b\eta$. We have that $\nu_{\mathcal C}(\omega_i^1)=\lambda_i$, since $\nu_{\mathcal C}(\mu x^ay^b\eta)\geq \lambda>\lambda_i$. In this way we produce an infinite list of strictly increasing divisorial order $1$-forms
 $$
 \omega_i=\omega_i^0,\omega_i^1,\omega_i^2,\ldots
 $$
such that $\nu_{\mathcal C}(\omega_i^j)=\lambda_i$, for any $i\geq 0$. For an index $j$ we have that $\nu_{E}(\omega_i^j)\geq c_\Gamma$ and then $\lambda_i\geq\nu_{E}(\omega_i^j)\geq c_\Gamma$ and this is a contradiction. So we necessarily have that $\nu_{\mathcal C}(\eta)=\lambda_i$. This ends the proof of Claim 3.
\vspace{5pt}

{\em Proof of Statements (1) and (2):} In view of Claim 1 and Claim 3, there is a $1$-form $\eta$ with $\nu_E(\eta)=t_i$ such that $\nu_\mathcal C(\eta)=\lambda_i$, whose initial part is proportional to the initial part of $x^{\ell_1}\omega_{i-1}$.
 In order to prove Statement (1), it remains to prove that if $\nu_E(\omega)=t_i$ then $\nu_{\mathcal C}(\omega)\leq \lambda_i$.
Assume that $\lambda=\nu_\mathcal C(\omega)>\lambda_i=\nu_\mathcal C(\eta)$. The $1$-form $\omega$ is basic and resonant and it has the same divisorial order as $\eta$. Hence there is a constant $\mu\ne 0$ such that
$$
\nu_E(\eta^1)>t_i=\nu_E(\eta)=\nu_E(\omega), \quad \eta^1=\eta-\mu\omega.
$$
The $1$-form $\eta^1$ satisfies that $\nu_\mathcal C(\eta^1)=\lambda_i\notin \Lambda_{i-1}$; by Claim 2, there are $a,b\geq 0$ and a constant $\mu'$ such that
$$
\nu_E(\eta^2)>\nu_E(\eta^1),\quad \eta^2=\eta^1-\mu 'x^ay^b\eta.
$$
We have that $\nu_{\mathcal C}(\eta^2)=\lambda_i$ and $\nu_E(\eta^2)>\nu_E(\eta^1)$. Repeating this procedure, we have a list of $1$-forms
 $
 \eta^1,\eta^2,\ldots
 $
 with strictly increasing divisorial order such that $\nu_{\mathcal C}(\eta^j)=\lambda_i$ for any $j$. We find a contradiction just by considering one of such $\eta^j$ with $\nu_E(\eta^j)\geq c_\Gamma$. This ends the proof of Statement (1).

 Let us prove Statement (2). Choose $\omega$ with $\nu_{\mathcal C}(\omega)=\lambda_i$.
 By Claim 2, we have that $\omega$ is reachable from $\eta$ and hence $\nu_E(\omega)\geq t_i$.   Assume by contradiction that $\nu_E(\omega)>t_i$. There is a constant $\mu$ and $a,b\geq 0$ with $a+b\geq 1$ such that
 $$
 \nu_E(\omega^1)>\nu_E(\omega),\quad \omega^1=\omega-\mu x^ay^b\eta.
 $$
Since
 $
 \nu_{\mathcal C}(\mu x^ay^b\eta)=an+bm+\lambda_i>\lambda_i
 $,
 we have that $\nu_\mathcal C(\omega^1)=\lambda_i$. Repeating the argument, we get a sequence of $1$-forms
 $
 \omega^0=\omega,\omega^1,\ldots
 $
with strictly increasing divisorial order such that $\nu_\mathcal C(\omega^j)=\lambda_i$ for any $j$. This is a contradiction.
\vspace{5pt}

 {\em Proof of  Statement (3): } By Claim 2, we have that $\omega_i$ is reachable from $x^{\ell_1}\omega_{i-1}$. By Statement (2) (already proved) we have that $\nu_E(\omega_i)=t_i$. Hence the initial part of $\omega_{i}$ is proportional to the initial part of $x^{\ell_1}\omega_{i-1}$. Consider a $1$-form $\omega$ with $\nu_{\mathcal C}(\omega)\notin \Lambda_{i-1}$. By Claim 2 the $1$-form $\omega$ is reachable from $x^{\ell_1}\omega_{i-1}$ and hence it is reachable from $\omega_i$. Then, there are $a,b\geq 0$ such that
 $$
 \nu_E(x^ay^b\omega_i)=an+bm+t_i=\nu_E(\omega).
 $$
 Since $\nu_\mathcal C(\omega)\notin \Lambda_{i-1}$, we have that $nm>\nu_\mathcal C(\omega)>\nu_E(\omega)>an+bm$, this implies the uniqueness of $a,b$.
  Let us show that
    $
    \nu_{\mathcal C}(\omega)\geq an+bm+\lambda_i
    $.
    Assume by contradiction that $
    \nu_{\mathcal C}(\omega)< an+bm+\lambda_i
    $. Consider
    $\omega^1=\omega-\mu x^ay^b\omega_{i}$
    such that $\nu_E(\omega^1)>\nu_E(\omega)$. In view of the contradiction hypothesis, we have that $\nu_{\mathcal C}(\omega^1)=\nu_{\mathcal C}( \omega)$. Moreover, if $\nu_E(\omega^1)=\nu_E(x^{a_1}y^{b_1}\omega_i)$ we also have that
    $
    \nu_{\mathcal C}(\omega)< a_1n+b_1m+\lambda_i
    $.
    The situation repeats and we obtain an infinite sequence of $1$-forms $\omega^0=\omega,\omega^1,\omega^2,\ldots$ with strictly increasing divisorial orders, such that $\nu_{\mathcal C}(\omega^j)=\nu_{\mathcal C}(\omega)$ for any $j\geq 0$. This is a contradiction.
\vspace{5pt}

{\em Proof of  Statement (4):} Noting that $\nu_E(x^{\ell_1}\omega_{i-1})=t_i$, by Statement (1) we have
$
\lambda_i\geq \nu_\mathcal C(x^{\ell_1}\omega_{i-1})=n\ell_1+\lambda_{i-1}=u_i
$. On the other hand, since $\lambda_i\notin \Lambda_{i-1}$, we have that $\lambda_i\ne u_i$ and hence $\lambda_i>u_i$.
\vspace{5pt}

{\em Proof of Statement (5):}
Consider $k=\lambda_i+na+mb$. Assume first that $k\notin \Lambda_{i-1}$. Let $\omega$ be such that $\nu_{\mathcal C}(\omega)=k$. We have to prove that $$\nu_{E}(\omega)\leq \nu_E(x^ay^b\omega_i)=an+bm+t_i.$$
    In view of Statement (3), we know that $\omega$ is reachable from $\omega_{i}$. Hence there are $a',b'\geq 0$ and a constant $\mu$ such that $\nu_E(\omega-\mu x^{a'}y^{b'}\omega_i)>\nu_E(\omega)$. Hence
    $$
    \nu_E(\omega)=\nu_E(x^{a'}y^{b'}\omega_i)=a'n+b'm+t_i.
    $$
    Assume by contradiction that
    $
    \nu_{E}(\omega)>\nu_E(x^ay^b\omega_i)=an+bm+t_i
    $.
    This implies that $a'n+b'm>an+bm$ and thus
    $$
    \nu_{\mathcal C}(x^{a'}y^{b'}\omega_i)=a'n+b'm+\lambda_i>k=an+bm+\lambda_i=\nu_{\mathcal C}(x^ay^b\omega_i)=\nu_{\mathcal C}(\omega).
    $$
   Put
    $\omega^1=\omega-\mu x^{a'}y^{b'}\omega_i$. We have that
    $
    \nu_{\mathcal C}(\omega^1)=k
    $.
    Repeating the argument with $\omega^1$, we obtain an infinite list of increasing divisorial orders $1$-forms
    $
    \omega^0=\omega,\omega^1,\omega^2,\ldots
    $
   such that $\nu_\mathcal C(\omega^j)=k\notin \Lambda_{i-1}$. This is a contradiction.

   Assume now that $k\in \Lambda_{i-1}$. There is an index $\ell\leq i-1$ such that
   $$
   k=an+bm+\lambda_{i}=a'n+b'm+\lambda_{\ell}.
   $$
   By Lemma \ref{lema:lambdamenoslambda},  we have that
   $
   \lambda_i-\lambda_\ell>t_i-t_\ell
   $
 and hence $ an+bm+t_{i}<a'n+b'm+t_{\ell}.$  The $1$-form $x^{a'}y^{b'}\omega_\ell$ satisfies that $k=\nu_\mathcal C(x^{a'}y^{b'}\omega_\ell)$ and
 $$
 \nu_E(x^{a'}y^{b'}\omega_\ell)=a'n+b'm+t_{\ell}>an+bm+t_{i}=\nu_E(x^ay^b\omega_i).
 $$
 This ends the proof.
 \end{proof}

\section{Delorme's decompositions}\label{apendice-Delorme-decomp}
In this Appendix, we provide a proof, using another approach, of Delorme's decompositions stated in Theorem \ref{th:delormedecomposition}. That is, we consider a cusp $\mathcal C\in \operatorname{Cusps}(E)$,
an  extended standard basis $\omega_{-1},\omega_0,\omega_1,\ldots,\omega_s;\omega_{s+1}$ of $\mathcal C$, where
$
\Lambda=\Gamma(n,m,\lambda_1,\ldots, \lambda_s)
$ is
the semimodule of differential values of $\mathcal C$. We have to prove that for any indices $0\leq j\leq i\leq s$, there is a decomposition
$$
\omega_{i+1}=\sum_{\ell=-1}^{j}f_\ell^{ij}\omega_\ell
$$
such that, for any $-1\leq \ell\leq j $ we have $\nu_{\mathcal C}(f_\ell^{ij}\omega_\ell)\geq v_{ij}$, where
$
v_{ij}=t_{i+1}-t_j+\lambda_j
$
and there is exactly one index $-1\leq k\leq j-1$ such that
$
\nu_{\mathcal C}(f_k^{ij}\omega_k)=\nu_{\mathcal C}(f_j^{ij}\omega_j)= v_{ij}
$.

Note that the case $i=0$ is straightforward. Indeed, we have $v_{00}=n+m$ and if we write $\omega_1=adx+bdy$, we necessarily have that $\nu_{\mathcal C}(adx)=\nu_\mathcal C(bdy)=n+m$ in view of the fact that the initial part of $\omega_1$ is proportional to $mydx-nxdy$.

Thus, we assume that $i\geq 1$.

\begin{lema}
\label{lema:clavedelorme}
  Given a $1$-form $\eta$ with $\nu_\mathcal C(\eta)>u_{i+1}$ and $\nu_E(\eta)>t_{i+1}$, we have
\begin{enumerate}
\item[a)] If
$\nu_E(\eta)<nm$, there is a $1$-form $\alpha$ such that $\nu_E(\eta-\alpha)>\nu_E(\eta)$ that can be decomposed as
$
\alpha=\sum_{\ell=-1}^{i}g_\ell\omega_\ell
$,
where  $\nu_\mathcal C(g_\ell\omega_\ell)>u_{i+1}$ for $-1\leq \ell\leq i$.
\item [b)] If $\nu_E(\eta)\geq nm$, there is a decomposition
$
\eta=\sum_{\ell=-1}^{i}h_\ell\omega_\ell
$
where each summand $h_\ell\omega_\ell$ satisfies that $\nu_\mathcal C(h_\ell\omega_\ell)>u_{i+1}$.
\end{enumerate}
\end{lema}
\begin{proof}
b) Assume that $\nu_E(\eta)\geq nm$, we have
$
\eta=fdx+gdy=f\omega_{-1}+g\omega_0
$,
where $\nu_E(fdx)\geq nm$ and $\nu_E(gdy)\geq nm$.  In view of Lemma
 \ref{lema:cotaconductor}, we have that $u_{i+1}<c_\Gamma+n<nm$.
We are done by taking the decomposition $
\eta=f\omega_{-1}+g\omega_0
$.

a) Assume that $\nu_E(\eta)<nm$. Note that $\eta$ is a basic $1$-form. There are two possible cases: $\eta$ is resonant or not.
Assume first that $\eta$ is not resonant. Then $$\nu_\mathcal C(\eta)=\nu_E(\eta)=\nu_E(\alpha)>u_{i+1},$$ where $\alpha$ is the initial part of $\eta$. Note that $\nu_E(\eta-\alpha)>\nu_E(\eta)$. We can write
$
\alpha=g_{-1}dx+g_0dy=g_{-1}\omega_{-1}+g_0\omega_0
$,
where
$$
\nu_\mathcal C(g_\ell\omega_\ell)\geq \nu_E(g_\ell\omega_\ell)\geq \nu_E(\alpha)=\nu_E(\eta)=\nu_\mathcal C(\eta)>u_{i+1},\quad  \ell=-1,0.
$$
This is the desired decomposition.

Assume that $\eta$ is resonant. Define
$
k=\max\{\ell\leq i; \; \eta \text{ is reachable from } \omega_\ell\}
$.
The fact that $\eta$ is resonant implies that $k\geq 1$ (recall that $i\geq 1$). Consider $a,b\geq 0$ and a constant $\varphi$ such that
\begin{equation}
\label{eq:omegak}
\nu_E(\tilde\eta)>\nu_E(\eta),\quad \tilde\eta=\eta-\varphi x^ay^b\omega_k.
\end{equation}
If we show that $\nu_\mathcal C(x^ay^b\omega_k)>u_{i+1}$, we are done. Let us do it. Assume first that $k=i$. We know that
$$
\nu_E(\eta)=\nu_E(x^ay^b\omega_{i})=an+bm+t_{i}>t_{i+1}=t_{i}+u_{i+1}-\lambda_{i}.
$$
This implies that $an+bm>u_{i+1}-\lambda_{i}$ and then $
\nu_{\mathcal C}(x^ay^b\omega_{i})=an+bm+\lambda_{i}>u_{i+1}
$.

Assume now that $1\leq k\leq i-1$. Let us reason by contradiction assuming that $\nu_{\mathcal C}(x^ay^b\omega_k)\leq u_{i+1}$. Denote by $\tilde\eta=\eta-\varphi x^ay^b\omega_k$. By Equation \eqref{eq:omegak}, we know that
\begin{equation}
\label{eq:nuetildeeta}
\nu_E(\tilde \eta)>\nu_E(x^ay^b\omega_k)=an+bm+t_k.
\end{equation}
Since $\nu_\mathcal C(x^ay^b\omega_k)\leq u_{i+1}<\nu_{\mathcal C}(\eta)$, we have that
\begin{equation}
\label{eq:nuctildeeta}
\nu_\mathcal C(\tilde \eta)=\nu_\mathcal C(x^ay^b\omega_k)=an+bm+\lambda_k.
\end{equation}
In view of Equations \eqref{eq:nuctildeeta} and \eqref{eq:nuetildeeta}, we can apply Statement (5) in Proposition \ref{prop:bs:delorme} to conclude that
\begin{equation}
an+bm+\lambda_k\in \Lambda_{k-1}.
\end{equation}
Let $\ell_1$ and $\ell_2$ be the $\Lambda_k$-limits.
Since $an+bm+\lambda_k\in \Lambda_{k-1}$, we have that $a\geq \ell_1$ or $b\geq \ell_2$,  by Lemma  \ref{lema:limits}.
We have four cases to be considered:
\begin{eqnarray*}
u_{k+1}= n\ell_1+\lambda_k \text{ and } a\geq \ell_1&;& u_{k+1}= n\ell_1+\lambda_k \text{ and } b\geq \ell_2;\\
u_{k+1}= m\ell_2+\lambda_k \text{ and } a\geq \ell_1&;& u_{k+1}= m\ell_2+\lambda_k \text{ and } b\geq \ell_2.
\end{eqnarray*}
Assume that $u_{k+1}= n\ell_1+\lambda_k$  and  $a\geq \ell_1$. This implies that
$x^ay^b\omega_k$, and hence $\eta$, is reachable from $x^{\ell_1}\omega_k$ and hence from $\omega_{k+1}$. This contradicts the maximality of the index $k$.

Assume that $u_{k+1}= n\ell_1+\lambda_k$  and  $b\geq \ell_2$ and $a<\ell_1$. We have that
$$
\nu_\mathcal C(x^ay^b\omega_k)\geq \nu_{\mathcal C}(y^{\ell_2}\omega_k)=m\ell_2+\lambda_k>u_{k+1}=n\ell_1+\lambda_k.
$$
Let $q_1$ and $q_2$ be the tops of $\Lambda_{k}$ and $Q_{\Lambda_k}$ the main top, we have that
$$
\nu_\mathcal C(x^ay^b\omega_k)\geq nQ_{\Lambda_k}\geq c_{\Lambda_k}+n.
$$
Recall that $c_{\Lambda_k}\leq n(Q_{\Lambda_k}-1)$ in view of Proposition \ref{lema:conductor1} . On the other hand, we know by Lemma  \ref{lema:cotaconductor} that
$
u_{i+1}<c_{\Lambda_i}+n\leq c_{\Lambda_k}+n
$.
We have the contradiction $u_{i+1}<c_{\Lambda_k}+n \leq \nu_\mathcal C(x^ay^b\omega_k)\leq u_{i+1}$.

The two remaining cases with $u_{k+1}= m\ell_2+\lambda_k$ may be considered in a similar way to the previous ones.
\end{proof}
\begin{proposition}
\label{prop:delorme} We can write
$
\omega_{i+1}=\sum_{\ell=-1}^{i}f_\ell\omega_\ell
$
where  $\nu_\mathcal C(f_\ell\omega_\ell)\geq u_{i+1}$ for $-1\leq \ell\leq i$ and such that $\nu_{\mathcal C}(f_{i}\omega_{i})=u_{i+1}$ and there is exactly one index $k\in \{-1,0,1,\ldots,i-1\}$ satisfying that $\nu_{\mathcal C}(f_k\omega_k)= u_{i+1}$.
\end{proposition}
\begin{proof} Let us consider first the case $i=0$. We know that $u_1=t_1=n+m$ and that $\omega_1$ is basic resonant, with $\nu_{E}(\omega_1)=n+m$. Then, there is a constant $\mu$ such that $\nu_E(\eta)>n+m$, where
$
\eta=\omega_1-\mu(mydx-nxdy)
$.
We can write
$
\eta=g_{-1}dx+g_0dy$, where $\nu_E(g_{-1})>m$ and  $\nu_E(g_0)>n$.
Let us put $f_{-1}= \mu my+g_{-1}$ and $f_0=-\mu nx+g_{0}$. We have that
$
\nu_{\mathcal C}(f_{-1})=m$ and $\nu_{\mathcal C}(f_0)=n
$;
hence
$
\omega_1=f_{-1}dx+f_0dy= f_{-1}\omega_{-1}+f_0\omega_0
$
is the desired decomposition.
\vspace{5pt}

Assume now that $1\leq i\leq s$. Let $\ell_1,\ell_2$ be the limits of $\Lambda_i$. By Remark \ref{rk:expresionesui}, there is exactly one index $k$ with $-1\leq k\leq i-1$ such that
\begin{enumerate}
\item  If $u_{i+1}=\ell_1n+\lambda_i$, then $u_{i+1}=\lambda_k+bm$ (note that $b\geq 1$).
\item  If $u_{i+1}=\ell_2m+\lambda_i$, then $u_{i+1}=\lambda_k+an$ (note that $a\geq 1$).
\end{enumerate}
Assume that $u_{i+1}=\ell_1n+\lambda_i$, the case $u_{i+1}=\ell_2m+\lambda_i$ is symmetric to this one.  We have that
\begin{equation}
\label{eq:valordiftiralui}
\nu_{\mathcal C}(x^{\ell_1}\omega_{i})=u_{i+1}=\nu_{\mathcal C}(y^{b}\omega_k).
\end{equation}
On the other hand, we have that
\begin{eqnarray*}
\nu_E(x^{\ell_1}\omega_{i})&=&t_i+ \ell_1n= t_{i}+(u_{i+1}-\lambda_{i})=t_{i+1};\\
\nu_E(y^b\omega_{k})&=&t_k+bm= t_k+u_{i+1}-\lambda_{k}.
\end{eqnarray*}
By Lemma \ref{lema:lambdamenoslambda}, we have that
$$
t_k+u_{i+1}-\lambda_{k}-t_{i+1}= (\lambda_{i}-\lambda_k)-(t_{i}-t_k)>0,
$$
and hence $t_{i+1}=\nu_E(x^{\ell_1}\omega_{i})<\nu_E(y^b\omega_{k})$.  Since both $\omega_{i+1}$ and $x^{\ell_1}\omega_{i}$ are basic resonant with the same divisorial order $t_{i+1}$, there is a constant $\varphi$ such that
$$\nu_E(\omega_{i+1}-\varphi x^{\ell_1}\omega_{i})>\nu_E(\omega_{i+1})=t_{i+1}.$$
 By Equation \eqref{eq:valordiftiralui},
there is a constant $\mu$ such that
$$
\nu_\mathcal C(\omega_{i+1}^{0})>u_{i+1}, \text{ where }  \omega_{i+1}^{0}=\varphi x^{\ell_1}\omega_{i}-\mu y^b\omega_{k}.
$$
Put $\eta^0=\omega_{i+1}-\omega_{i+1}^0=\omega_{i+1}-\varphi x^{\ell_1}\omega_{i}+\mu y^b\omega_{k}$. We have that
$
\nu_E(\eta^0)>t_{i+1}$ and $\nu_{\mathcal C}(\eta^0)>u_{i+1}
$
in view of the following facts:
\begin{enumerate}
\item $\nu_E(\eta^0)\geq\min\{\nu_E(\omega_{i+1}-\varphi x^{\ell_1}\omega_{i}),\nu_E(\mu y^b\omega_{k})\}>t_{i+1}$.
\item $
\nu_{\mathcal C}(\eta^0)\geq \min \{\nu_\mathcal C(\omega_{i+1}),\nu_{\mathcal C}(\omega_{i+1}^0)\}= \min \{\lambda_{i+1},\nu_{\mathcal C}(\omega_i^0)\}>u_{i+1}$. Recall that $\Lambda$ is an increasing semimodule; (here we put $\lambda_{s+1}=\infty$).
\end{enumerate}
 The proof is now a consequence of Lemma \ref{lema:clavedelorme} as follows.
 We start with $\eta^0$ as before.
If $\nu_E(\eta^0)\geq nm$, we apply Lemma \ref{lema:clavedelorme} b).  We are done by taking the decomposition
$$
\omega_{i+1}=\omega_{i+1}^0+ \sum_{\ell=-1}^{i}h_\ell\omega_\ell= \varphi x^{\ell_1}\omega_{i}-\mu y^b\omega_{k}+ \sum_{\ell=-1}^{i}h_\ell\omega_\ell.
$$
If $\nu_E(\eta^0)<nm$, we apply Lemma \ref{lema:clavedelorme} a) and we obtain $\eta^1=\eta^0-\sum_{\ell=-1}^{i}g_\ell\omega_\ell$ such that $\nu_E(\eta^1)>\nu_E(\eta^0)>t_{i+1}$ and $\nu_\mathcal C(\eta^1)>u_{i+1}$. If $\nu_E(\eta^1)< nm$, we re-apply Lemma \ref{lema:clavedelorme} a) to $\eta^1$. After finitely many steps, we obtain
$$
\tilde\eta=\eta^0- \sum_{\ell=-1}^{i}\tilde g_\ell\omega_\ell,
$$
such that $\nu_E(\tilde \eta)\geq nm$, $\nu_{\mathcal C}(\tilde\eta)>u_{i+1}$ and $\nu_{E}(\tilde\eta)>t_{i+1}$ ,
where $\nu_{\mathcal C}(\tilde g_\ell\omega_\ell)>u_{i+1}$ for any $-1\leq \ell\leq i$. We apply Lemma \ref{lema:clavedelorme} b) to $\tilde\eta$ to obtain that $\tilde\eta=\sum_{\ell=-1}^{i}\tilde h_\ell\omega_\ell$ with $\nu_{\mathcal C}(\tilde h_\ell\omega_\ell)>u_{i+1}$. The desired decomposition is given by
$$
\omega_{i+1}= \omega_{i+1}^0 + \sum_{\ell=-1}^{i}(\tilde g_\ell+\tilde h_\ell)\omega_\ell= \varphi x^{\ell_1}\omega_{i}-\mu y^b\omega_{k} + \sum_{\ell=-1}^{i}(\tilde g_\ell+\tilde h_\ell)\omega_\ell.
$$
This ends the proof.
\end{proof}

Let us end the proof of Theorem \ref{th:delormedecomposition}. We already know that it is true when $j=i$, in view of Proposition \ref{prop:delorme}. We assume that the result is true for $j+1\leq i$ and let us show that it is true for $0\leq j\leq i$. In order to simplify notations, let us write
$$
v_j= t_{i+1}-t_j+\lambda_j ,\quad v_{j+1}= t_{i+1}-t_{j+1}+\lambda_{j+1}.
$$
Recall that $v_{j+1}=v_j+\lambda_{j+1}-u_{j+1}$.
By induction hypothesis, we have that
$$
\omega_{i+1}=\sum_{\ell=-1}^{j+1}h_\ell\omega_{\ell},
$$
where $\nu_\mathcal C(h_\ell\omega_{\ell})\geq v_{j+1}$ for any $-1\leq \ell\leq j+1$ and
$
\nu_\mathcal C(h_{j+1}\omega_{j+1})= v_{j+1}
$.
We apply Proposition \ref{prop:delorme} to write
$
\omega_{j+1}=\sum_{\ell=-1}^{j}g_\ell\omega_\ell
$,
where $\nu_\mathcal C(g_\ell\omega_\ell)\geq u_{j+1}$ for any $-1\leq \ell\leq j$ and there is exactly one index $k$ such that
$
\nu_\mathcal C(g_{j}\omega_{j})=\nu_\mathcal C(g_{k}\omega_{k})= u_{j+1}
$.
Now, we have an expression
$$
\omega_{i+1}=\sum_{\ell=-1}^{j}f_\ell\omega_\ell,\quad f_\ell=h_\ell+h_{j+1}g_\ell.
$$
We have the following properties:
\begin{enumerate}
\item $ \nu_\mathcal C(h_\ell\omega_\ell)> v_{j}$, for any $-1\leq \ell\leq j$. Indeed, we know that
    $$
    v_{j+1}= v_j+(\lambda_{j+1}-u_{j+1})>v_j,
    $$
    recall that the semimodule is increasing and then $\lambda_{j+1}>u_{j+1}$.
\item  $\nu_\mathcal C(h_{j+1}g_\ell\omega_\ell)\geq v_{j}$ and $k, j$ are the unique indices such that
    $$
    \nu_\mathcal C(h_{j+1}g_j\omega_j)=\nu_\mathcal C(h_{j+1}g_k\omega_k)=v_j.
    $$
    In order to prove this, it is enough to note that
    $$
\nu_\mathcal C(h_{j+1}g_\ell\omega_\ell)=(v_{j+1}-\lambda_{j+1})+\nu_\mathcal C(g_\ell\omega_\ell)\geq (v_{j+1}-\lambda_{j+1})+u_{j+1}=v_j
$$
and the equality holds exactly for the indices $\ell=j,k$.
\end{enumerate}
 The desired result comes from the above properties (1) and (2), noting that
 $$
 \nu_\mathcal C(f_\ell\omega_\ell)\geq \min\{\nu_\mathcal C(h_\ell\omega_\ell),
 \nu_\mathcal C(h_{j+1}g_\ell\omega_\ell)\}
 $$
and the equality holds when the two values are different.

\end{document}